\documentclass[1p]{elsarticle}
\biboptions{numbers,sort&compress}

\usepackage{xcolor}

\usepackage{subfig}
\usepackage{float}
\usepackage{amssymb}
\usepackage{amsmath} 

\begin{document}

\title{A $p$-Multigrid Accelerated Nodal Spectral Element Method for Free-Surface Incompressible Navier-Stokes Model of Nonlinear Water Waves}


\author[1]{Anders Melander\corref{cor1}}
\ead{adame@dtu.dk}

\author[2]{Wojciech Laskowski}
\ead{lskwj@gmail.com}

\author[3]{Spencer J. Sherwin}
\ead{s.sherwin@imperial.ac.uk}

\author[1]{Allan P. Engsig-Karup}
\ead{apek@dtu.dk}

\cortext[cor1]{Corresponding author}
\affiliation[1]{organization={Department of Applied Mathematics and Computer Science, Center for Energy Resources Engineering (CERE), Technical University of Denmark}, addressline={Richard Petersens Plads}, city={2800 Kgs. Lyngby}, country={Denmark}}
\affiliation[2]{organization={ETSIAE-UPM (School of Aeronautics - Universidad Politécnica de Madrid)}, addressline={Plaza de Cardenal Cisneros 3}, cite={28040 Madrid}, country={Spain}}
\affiliation[3]{organization={Department of Aeronautics, Imperial College London}, addressline={Exhibition Road}, city={SW7 2AZ London}, country={England}}


\begin{abstract}
We present a spectral element model for general-purpose simulation of non-overturning nonlinear water waves using the incompressible Navier-Stokes equations (INSE) with a free surface. The numerical implementation of the spectral element method is inspired by the related work by Engsig-Karup et al. \cite{Allan16} and is based on nodal Lagrange basis functions, mass matrix-based integration and gradient recovery using global $L^2$ projections. The resulting model leverages the high-order accurate 
 -- possibly exponential -- error convergence and has support for geometric flexibility allowing for computationally efficient simulations of nonlinear wave propagation. An explicit fourth-order accurate Runge-Kutta scheme is employed for the temporal integration, and a mixed-stage numerical discretization is the basis for a pressure-velocity coupling that makes it possible to maintain high-order accuracy in both the temporal and spatial discretizations while preserving mass conservation. Furthermore, the numerical scheme is accelerated by solving the discrete Poisson problem using an iterative solver strategy based on a geometric $p$-multigrid method. This problem constitutes the main computational bottleneck in INSE models. It is shown through numerical experiments, that the model achieves spectral convergence in the velocity fields for highly nonlinear waves, and there is excellent agreement with experimental data for the simulation of the classical benchmark of harmonic wave generation over a submerged bar. The geometric $p$-multigrid solver demonstrates $O(n)$ computational scalability simulations, making it a suitable efficient solver strategy as a candidate for extensions to more complex, real-world scenarios.
\end{abstract}


\maketitle
\section{Introduction}
The ability to accurately model nonlinear water waves is an important aspect in many fields such as oceanography, coastal and offshore engineering, and environmental science. The ability to use numerical simulations of unsteady water waves and flow kinematics is therefore of large interest in these fields. Such simulations can allow for the efficient study of flow characteristics such as wave loads. For simulating unsteady water waves, the most widely used models are based on Navier-Stokes equations (NSE), or on fully nonlinear potential flow (FNPF) models. While FNPF models can be significantly more efficient than full CFD solvers based on NSE \cite{Ransley2019}, they do not take into account viscous and rotational effects and have a tendency to over-predict the wave heights in wave-structure interaction applications due to lack of viscous effects \cite{SriramEtAl2021}. Hence, in such cases the use of models based on the incompressible Navier-Stokes equations are required and efficient numerical schemes are needed.

The full set of NSE can been solved in many ways using different numerical schemes. Some of the conventional methods are based on implicit surface-capturing techniques such as level-set techniques \cite{Osher88,GroossHesthaven06,BIHS2016191} and the volume of fluid (VOF) method \cite{Hirt81} that has been popularized through open-source frameworks such as OpenFoam \cite{WellEtAl1998}. While these methods have the benefit of being able to handle breaking waves, they are disadvantaged in several ways that affect the computational cost and overall accuracy of the numerical schemes. Despite being widely used, the methods incur high computational costs which limits their use in large scale CFD simulations in both space and time. Moreover, applying pressure boundary conditions to the free surface is inaccurate due to the boundary crossing computational elements \cite{Lin02}.  Alternatively, another widely used modeling method is to track the free surface explicitly through a kinematic boundary condition. This way, the free surface elevation is modeled, which opens for designing efficient computational efficient schemes, however, at the price of not being able to handle breaking waves. By explicit tracking the free surface by incorporating a kinematic free surface boundary condition it is possible to keep track of the shape of the spatial domain at all times. Common examples of this is the arbitrary Lagrangian-Eulerian (ALE) method \cite{Donea04}, the mixed Eulerian-Lagrangian method \cite{LonguetHigginsCokelet1976,EngsigKarupMontesering2019} or the $\sigma$-transform method used to transform the vertical coordinate to a regular computational domain, which was first used for CFD wave simulation by Li \& Flemming \cite{Li01}. It is noted that Decoene \& Gerbeau \cite{Decoene09} showed that the $\sigma$-transform can be seen as a special case of the ALE method. Li \& Flemming used a 2nd order accurate numerical scheme in both time and space, and coupled pressure and velocity through a projection method, first suggested by Chorin \cite{Chorin68}. Another popular pressure-velocity coupling technique is the Harlow \& Welch algorithm \cite{HarlowWelch1965} that uses a predictor-corrector method to compute velocity and pressure fields, ensuring that the free surface satisfies kinematic and dynamic conditions while maintaining incompressibility.

Recently, higher-order numerical methods have also been used for the $\sigma$-transformed problem, such as Engsig-Karup et al \cite{Allan24} who used a high-order finite difference scheme to solve the free surface Navier-Stokes problem in 2D using a high-order time-stepping method and resolving the pressure through a mixed-stage Poisson problem to achieve a divergence free numerical scheme. Moreover, Pan et al \cite{Pan21} solved the NSE using a discontinuous Galerkin finite element method (DGFEM), while Melander et al. \cite{Melander24} designed a pseudo-spectral method, both achieving spectral convergence of the error. Spectral methods have been shown to also be computationally efficient for other type of wave models such as Boussinesq models \cite{MADSEN_BINGHAM_LIU_2002,Allan06} and FNPF models \cite{Christiansen13,Klahn20,Bonnefoy10,MelanderEngsigKarup2024}. Patera \cite{patera84} developed the spectral element method (SEM) for fluid dynamics, combining the attractive error convergence of spectral methods, with the geometric flexibility of the finite element method, cf. the original review of using SEM for INS \cite{MadayPatera1989}. Robertson and Sherwin \cite{sherwin99} presented a SEM discretization of the FNPF model that shows stability issues related to certain mesh configurations. This issue was mitigated in the work of Engsig-Karup et al \cite{Allan16} and stabilised free surface simulations using SEM for the FNPF model was presented. The use of high-order numerical schemes allows for high accuracy and low numerical diffusion at a lower cost for a given accuracy compared to more traditional low-order methods \cite{Kreiss72}. Hence, such schemes are very attractive choices for time-dependent problems where numerical diffusion errors are important to minimize due to build-up of errors in the settings of both large-scale simulations and longer time integration. In relation to these previous works, a motivation behind this work is to apply the SEM to develop a new numerical free surface NSE-based model that can form a basis for wave propagation and wave-structure interaction applications. 

In CFD solvers, the Poisson problem for pressure often represents a significant computational bottleneck, critically affecting runtime performance. As a result, efficient iterative solvers are essential for effectively resolving the linear systems that arise from the numerical discretization. Iterative multigrid solvers, for instance, are among the most efficient and $\mathcal{O}(n)$-scalable algorithms. These solvers capitalize on the ability of traditional stationary iterative methods to efficiently reduce high-frequency errors \cite{Trottenberg01} and are widely used for accelerating advanced CFD solvers \cite{brandt1998barriers,JamesEtAl2002}. By employing a solution approach that solves systems of algebraic equations across multiple grid levels, the effectiveness of these methods is preserved, as low-frequency errors on finer grids become high-frequency errors on coarser grids. Geometric multigrid methods have demonstrated their effectiveness in water wave problems, such as free surface NSE modeling \cite{Li01} and FNPF modeling \cite{Allan09}. Furthermore, the use of high-order accurate numerical methods facilitates the application of geometric $p$-multigrid methods, extending multigrid techniques to leverage a hierarchy of varying polynomial orders. Engsig-Karup and Laskowski \cite{AllanLaskowski} demonstrate efficient wave propagation and wave-structure interaction by using a geometric $p$-multigrid method as a preconditioner for the stationary defect correction (PDC) method and preconditioned conjugate gradient (PCG) methods in FNPF models discretized with SEM. Moreover, Melander et. al \cite{Melander24} showed that a $p$-multigrid method could achieve improved computational efficiency when used to solve the Poisson problem in a free Surface Navier-Stokes model discretized with a pseudo-spectral numerical method. 

\subsection{Paper contributions}
In this work, we present a novel free surface incompressible Navier-Stokes solver based on the high-order nodal spectral element method. The work can be seen as a multi-domain extension of the single-domain pseudospectral method described in \cite{Melander24}. The design of a high-order numerical scheme is a basis for cost-efficient simulations that can achieve high accuracy and low numerical diffusion, along with support for geometric flexibility in comparison to the single domain approach. Through numerical experiments spectral (faster than algebraic fixed order) convergence is demonstrated together with longer time nonlinear wave propagation over uneven bathymetry. Moreover, to address the computational bottleneck of solving the Poisson problem, we propose to accelerate the iterative solver strategy with a geometric $p$-multigrid method. We show that this iterative solver strategy for the Poisson problem can achieve $O(n)$ computational scalability similar to the $p$-multigrid accelerated SEM scheme for the Laplace problem in the FNPF model due to \cite{Allan16, AllanLaskowski}, contributing to enabling cost-efficient and scalable solution of the Poisson problem. 

\subsection{Paper organization}
The paper has been organized as follows. In Section \ref{sec:goveq} we introduce the governing equations for the incompressible Navier-stokes flow with a free surface. In Section \ref{sec:sigma} we introduce the $\sigma$-transform for the coordinates. In Section \ref{sec:temporal} we detail the temporal discretization along with the velocity-pressure coupling achieved by mass conservation through a mixed-stage Poisson problem. In Section \ref{sec:spatial} we detail the spatial discretization using the Spectral Element Method. The weak formulations of the governing equations are presented along with element and basis construction. In Section \ref{sec:mg} we go through the $p$-multigrid method used for the occurring mixed-stage Poisson problem. Lastly, in Section \ref{sec:res} we presents results of the proposed solver both in terms of accuracy of simulations and computational efficiency of the geometric $p$-multigrid method.

\section{Governing equations}\label{sec:goveq}
\begin{figure}[H]
\vspace*{0cm}
\hspace*{-2.5cm}
\centering
  \includegraphics[width=160mm]{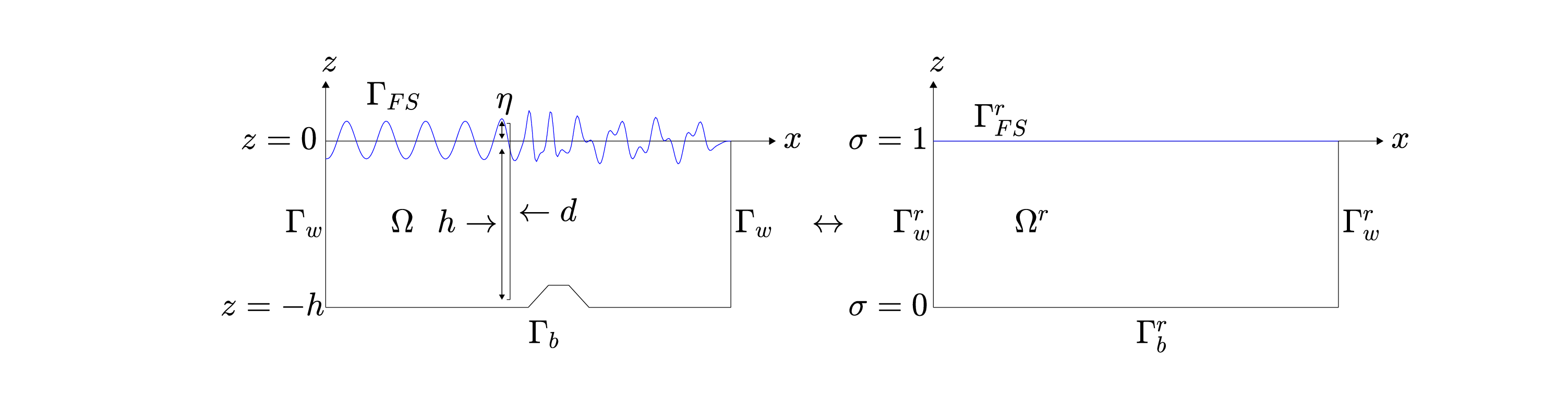}
  \caption{Illustration of the physical and $\sigma$-transformed domains.}
  \label{fig:setup1}
\end{figure}
The description of the evolution of water waves in the time domain can be described by the incompressible Navier-Stokes equations (INSE) along with a free surface kinematic boundary condition \cite{EGNL13,Allan24}. By assuming that the fluid density is constant, the fluid flow becomes divergence free. In the following, we introduce the fluid domain $\Omega\subset \mathbb{R}^d$ ($d=2$) be a bounded, connected domain with a piecewise smooth spatial domain boundary $\Gamma=\partial\Omega$. We introduce restrictions of the this boundary to the free surface $\Gamma_{FS}\subset \mathbb{R}^{d-1}$, solid impermeable domain boundaries such as walls $\Gamma_w$, and the bathymetry $\Gamma_b\subset\mathbb{R}^{d-1}$. The time domain is defined by $T:0\leq t \leq t_f$, where $t$ denotes the time variable and $t_f$ the final time of  simulation. Figure \ref{fig:setup1} shows an example of the domain, along with the $\sigma$-transformed domain to be introduced in Section \ref{sec:sigma}.

The INSE can be stated in terms of the mass conservation equation and momentum equations, and in the following the governing equations are formulated in terms of two spatial dimensions ($d=2$) and time in an Eulerian frame of reference. Find ${\bf u}$, $p$ such that  
\begin{subequations}
\begin{align}
    &\nabla\cdot\textbf{u} = 0, & \text{in} \;\;\; \Omega , 
    \label{eq:div1}\\
    &\frac{\partial\textbf{u}}{\partial t} =-\frac{1}{\rho}\nabla p  + \textbf{g} + \nu \nabla^2 \textbf{u}-\textbf{u}\cdot\nabla\textbf{u},& \text{in} \;\;\; \Omega \times T .\label{eq:mom1}
\end{align} 
Here $\textbf{u}=(u,w)^T$ are the velocities $[\frac{m}{s}]$ and $p$ is the pressure $[\frac{kg}{m^3}]$. Moreover, $\rho$ is the density of the fluid and assumed defined as $\rho=999.70$ $[\frac{kg}{m^3}]$ (at 10$^{\circ}$ C), $\nu$ is the kinematic viscosity $[\frac{m^2}{s}]$, and $\textbf{g}=(0,-g_z)$ is the gravitational acceleration where $g_z=9.81$ $[\frac{m^2}{s}]$ is assumed.\\
The free water surface is described by a kinematic boundary condition given as
\begin{align}
    \frac{\partial\eta}{\partial t} + \tilde{u}\frac{\partial\eta}{\partial x}  = \tilde{w} \quad \quad \quad \quad \quad \quad  & \text{on} \;\;\; \Gamma_{FS}, \label{eq:freesurf1}
\end{align}
\end{subequations}
where $\eta$ is the free surface, and '$\sim$' denotes a variable evaluated at the free surface level corresponding to $z=\eta({\bf x},t)$, e.g. the free surface velocities $\tilde{{\bf u}} = (\tilde{u},\tilde{w})^T = {\bf u}({\bf x,z=\eta},t)$. 

To solve the governing equations, we need to define spatial domain boundary conditions and relate velocity and pressure through a pressure-velocity coupling that satisfies mass conservation \eqref{eq:div1}.

\subsection{Boundary conditions}

For the solution of governing equations suitable boundary conditions are needed. Therefore, in the setting of a numerical wave tank, the following boundary conditions are defined.

We restrict the fluid to be contained within some domain $\Omega$, bounded from above by the kinematic free surface boundary condition (\ref{eq:freesurf1}). Moreover, the spatial domain is enclosed by wall boundaries and a bathymetry. At these boundaries, we assume an impermeability condition in the form of a slip condition, i.e.
\begin{align}    
\textbf{n}\cdot\textbf{u}=0\quad\text{on}\quad\Gamma_{w}\cup\Gamma_b.
\end{align}
When solving for the dynamic pressure via a pressure-velocity coupling scheme to be defined in section \ref{sec:poissonEq}, we define a reference level for the pressure $p$ at the free surface water line to be 
\begin{align}
    p=0\quad\text{on}\quad\Gamma_{FS}.
\end{align}
It is possible to split the pressure into its static and dynamic parts, such that
\begin{align}
    p = p_S+p_D,\quad p_S = \rho g(\eta-z). \label{eq:pressplit}
\end{align}
Remark, the hydrostatic pressure $p_S$ contribution is per definition also zero at the surface, meaning
\begin{align}
    p_S=0\quad\text{on}\quad\Gamma_{FS},\quad
    p_D=0\quad\text{on}\quad\Gamma_{FS}.
\end{align}
This splitting can be used in the derivation of the pressure-velocity coupling to solve for the unknown dynamic pressure contribution at an instant in time.

\subsection{The $\sigma$-coordinate transform}\label{sec:sigma}
The free surface evolves with time implying that the spatial domain for the fluid is time-dependent. This means that a possible re-meshing and update of computational operators has to happen in every time step. To avoid expensive re-meshing and reduce the cost of updating the spatial operations,  it is common to employ a $\sigma$-transform which maps the physical domain to a time-invariant reference domain $\Omega^r$. Since the physical spatial $\Omega$ only changes with respect to time in the vertical direction, the vertical mapping is defined as 
\begin{equation}
    t^* = t, \quad x^* = x,  \quad \sigma = \sigma(t,x,z). 
\end{equation}
Here '$*$' denotes the reference domain. In this domain, the vertical coordinate direction is given by $\sigma$ which ensures that the reference domain is time-invariant, 
\begin{equation}
    \sigma = \frac{z+h(x)}{d(x,t)}, \quad 0 \leq \sigma \leq 1. \label{eq:sigtrans}
\end{equation}
Here $h(x)$ denotes the still water height and $d(x,t)$ denotes the total water height, given as $d(x,t) = \eta(x,t) + h(x)$. 

The chain rule is employed to transform operations in the physical domain to corresponding operations in the reference domain. We know that for any function $f$, the following must hold true for the transformation to be consistent,
\begin{equation}
    f(t,x,z) = f(t^*,x^*,\sigma(t,x,z) ).
\end{equation}
We apply the chain rule for each of the four Cartesian variables on $f$ to find the transformed gradient operators 
    \begin{align}
        &\frac{\partial f}{\partial t} = \frac{\partial f}{\partial t^*} + \frac{\partial f}{\partial \sigma} \frac{\partial \sigma}{\partial t},\quad 
        \frac{\partial f}{\partial x} = \frac{\partial f}{\partial x^*} + \frac{\partial f}{\partial \sigma} \frac{\partial \sigma}{\partial x},\quad
        \frac{\partial f}{\partial z} = \frac{\partial f}{\partial \sigma} \frac{\partial \sigma}{\partial z}.
    \end{align}
We apply these relations to derive a $\sigma$-transformed version of the Navier-Stokes equations  (\ref{eq:div1})-(\ref{eq:mom1}) and introduce the pressure splitting (\ref{eq:pressplit}) to make the dependence on the dynamic pressure explicit in the formulation. Following \cite{Li01} the following variable is introduced responsible for accounting for the temporal changes of the spatial domain 
\begin{align}
    w_{\sigma} = \frac{\partial \sigma}{\partial t} + u\frac{\partial \sigma}{\partial x}  + w\frac{\partial \sigma}{\partial z},
\end{align}
 giving \begin{subequations} 
    \begin{align}
            & \frac{\partial u}{\partial x^*} + \frac{\partial u}{\partial \sigma}\frac{\partial \sigma}{\partial x} + \frac{\partial w}{\partial \sigma}\frac{\partial \sigma}{\partial z} = 0,\\ 
            & \frac{\partial u}{\partial t^*} + u\frac{\partial u}{\partial x^*} + w_{\sigma}\frac{\partial u}{\partial \sigma} = -\frac{1}{\rho}\left(\frac{\partial p_D}{\partial x^*}+\frac{\partial p_D}{\partial \sigma}\frac{\partial \sigma}{\partial x}+ \rho g\frac{\partial \eta}{\partial x^*}\right)\nonumber\\
            &+\nu\left(\frac{\partial^2 u}{\partial x^{*2} }+\frac{\partial^2 u}{\partial \sigma^2}\left(\frac{\partial \sigma}{\partial x}\right)^2+2\frac{\partial^2 u}{\partial x^* \partial \sigma}\frac{\partial \sigma}{\partial x}+\left(\frac{\partial^2 \sigma}{\partial x^{2}}  \right) \frac{\partial u}{\partial \sigma}+\frac{\partial^2 u}{\partial \sigma^2}\left(\frac{\partial \sigma}{\partial z}\right)^2\right),\\
            & \frac{\partial w}{\partial t^*} + u\frac{\partial w}{\partial x^*} +  w_{\sigma}\frac{\partial w}{\partial \sigma} = -\frac{1}{\rho}\left(\frac{\partial p_D}{\partial \sigma}\frac{\partial \sigma}{\partial z}\right)-g_z\nonumber\\
            &+\nu\left(\frac{\partial^2 w}{\partial x^{*2} }+\frac{\partial^2 w}{\partial \sigma^2}\left(\frac{\partial \sigma}{\partial x}\right)^2+2\frac{\partial^2 w}{\partial x^* \partial \sigma}\frac{\partial \sigma}{\partial x}+\left(\frac{\partial^2 \sigma}{\partial x^{2}}  \right) \frac{\partial w}{\partial \sigma}+\frac{\partial^2 w}{\partial \sigma^2}\left(\frac{\partial \sigma}{\partial z}\right)^2\right).
    \end{align}
\end{subequations}
Lastly, we introduce the modified gradient, Laplace operator and $\textbf{u}_{\sigma}$ vector defined as
\begin{subequations}
    \begin{align}
            &\nabla_{\sigma} = \left(\frac{\partial}{\partial x^*} + \frac{\partial \sigma}{\partial x}\frac{\partial}{\partial \sigma}, \frac{\partial \sigma}{\partial z}\frac{\partial}{\partial \sigma}\right)^T, \\
            & \nabla_{\sigma}^2 = \frac{\partial^2}{\partial x^{*2}}  + \left( \left(\frac{\partial \sigma}{\partial x}\right)^2  + \left(\frac{\partial \sigma}{\partial z}\right)^2 \right)\frac{\partial^2 }{\partial \sigma^2} \nonumber\\
            &+ 2 \frac{\partial \sigma}{\partial x} \frac{\partial^2}{\partial \sigma \partial x^*} + \left(\frac{\partial^2 \sigma}{\partial x^{2}}  \right) \frac{\partial }{\partial \sigma},\\
            &\textbf{u}_{\sigma} = (u,w_{\sigma})^T.
    \end{align}
\end{subequations}
Now, the $\sigma$-transformed Navier-Stokes equations are expressed on vector form similar to the original equations
\begin{subequations}
\begin{align}
    & \nabla_{\sigma} \cdot \textbf{u} = 0, \\
    &\frac{\partial \textbf{u}}{\partial t^*} + \textbf{u}_{\sigma}\cdot  \hat{\nabla}_\sigma \textbf{u} = - \frac{1}{\rho} (\nabla_{\sigma} p_D + \hat{\nabla}_\sigma p_S) + \textbf{g} + \nu \nabla_{\sigma}^2 \textbf{u}.\label{eq:sigmom1} 
\end{align}
\end{subequations}
Note that we have introduced the operator
\begin{align}
     &\hat{\nabla}_\sigma = \left(\frac{\partial}{\partial x^*},\frac{\partial}{\partial \sigma}\right)^T.
\end{align}
We evaluate the $\sigma$-dependent derivatives by applying the chain rule to the transformation function (\ref{eq:sigtrans}) and find
\begin{equation} \label{eq:1_Scoef}
    \begin{aligned}
        &\frac{\partial \sigma}{\partial t} = -d^{-1}\left(\sigma \frac{\partial d}{\partial t}\right),\quad 
        &\frac{\partial \sigma}{\partial x} = d^{-1}\left(\frac{\partial h}{\partial x}  -\sigma \frac{\partial d}{\partial x}  \right),\\
        &\frac{\partial^2 \sigma}{\partial x^2} = d^{-1}\left(\frac{\partial^2 h}{\partial x^2} - \sigma \frac{\partial^2 d}{\partial x^2}  - 2 \frac{\partial \sigma}{\partial x} \frac{\partial d}{\partial x} \right),  
        &\frac{\partial \sigma}{\partial z} =  d^{-1}. 
    \end{aligned} 
\end{equation}
It is noted that the physical derivatives that appear in the governing equations can be evaluated using the formulas presented in this section, and with the spatial gradients of the solution variables and free surface variables computed numerically using the spectral element method. Also, it is noted that  $\frac{\partial d}{\partial t}$ in (\ref{eq:1_Scoef}) is related directly to the kinematic boundary condition through
\begin{equation}
    \frac{\partial d\vert_{z=\eta}}{\partial t} = \frac{\partial \eta}{\partial t}+\frac{\partial h}{\partial t}=w\vert_{z=\eta} - u\vert_{z=\eta} \frac{\partial \eta}{\partial x} + \frac{\partial h}{\partial t}.
\end{equation} 
\subsection{The relationship between the $\sigma$-transform and the Arbitrary Lagrangian-Eulerian method}
To clarify that the $\sigma$-transform technique is appropriately handling the time derivative in the moving frame of reference, we show that the $\sigma$-transform can be considered a special case of the Arbitrary Lagrangian–Eulerian (ALE) method \cite{Donea04} for free surface flow. The ALE formulation of the free surface Navier-Stokes equations is done by transforming the time derivative to the ALE reference domain \cite{Duarte04}. We let $\hat{t}$, $\hat{\textbf{x}}$ refer to the variables in the ALE reference domain, and let $\textbf{x}=f(\hat{\textbf{x}})$ be the transformation. The time derivative is transformed as
\begin{align}
     \frac{\partial \textbf{u}}{\partial t}= \frac{\partial \textbf{u}}{\partial \hat{t}}-\frac{\partial \textbf{x}}{\partial t}\frac{\partial \textbf{u}}{\partial \textbf{x}},\label{eq:ALEdt2}
\end{align}
leading to the traditional ALE form of the momentum equations:
\begin{align}
    \frac{\partial\textbf{u}}{\partial \hat{t}}+(\textbf{u}-\textbf{v}_f)\cdot\nabla\textbf{u} =-\frac{1}{\rho}\nabla p  + \textbf{g} + \nu \nabla^2 \textbf{u}.\label{eq:ALEmom2}
\end{align}
Here $\textbf{v}_f = \frac{\partial \textbf{x}}{\partial t} = \frac{\partial }{\partial t}f(\hat{\textbf{x}})$ denotes the domain velocity. If instead, we transform the time derivative by use of the inverse transformation $\hat{\textbf{x}} = g(\textbf{x})$, the time derivative can be transformed as
\begin{align}
     \frac{\partial \textbf{u}}{\partial t}= \frac{\partial \textbf{u}}{\partial \hat{t}}+\frac{\partial \hat{\textbf{x}}}{\partial t}\frac{\partial \textbf{u}}{\partial \hat{\textbf{x}}}.\label{eq:ALEdt3}
\end{align}
This gives us the alternative ALE form:
\begin{align}
   \frac{\partial \textbf{u}}{\partial \hat{t}}+\frac{\partial \hat{\textbf{x}}}{\partial t}\frac{\partial \textbf{u}}{\partial \hat{\textbf{x}}}+ \textbf{u}\cdot\nabla\textbf{u} =-\frac{1}{\rho}\nabla p  + \textbf{g} + \nu \nabla^2 \textbf{u}.\label{eq:ALEmom3}
\end{align}
It is clear, that if only vertical mesh movement in the ALE frame is considered, the transformation can be described through the $\sigma$-transform, i.e. 
\begin{align}
    g(\textbf{x}) =\begin{pmatrix}
\hat{x}\\
\sigma(x,z,t)
\end{pmatrix}=\begin{pmatrix}
\hat{x}\\
\frac{z+h(x)}{d(x,t)}
\end{pmatrix}.
\end{align}
This gives us the time derivative in the $\sigma$-domain as 
\begin{align}
     \frac{\partial \textbf{u}}{\partial t}= \frac{\partial \textbf{u}}{\partial \hat{t}}+\frac{\partial \sigma}{\partial t}\frac{\partial \textbf{u}}{\partial \sigma},\label{eq:ALEsigdt}
\end{align}
and the ALE-$\sigma$ form corresponding to  (\ref{eq:ALEmom3}) as 
\begin{align}
   \frac{\partial \textbf{u}}{\partial \hat{t}}+\frac{\partial \sigma}{\partial t}\frac{\partial \textbf{u}}{\partial \sigma}+ \textbf{u}\cdot\nabla\textbf{u} =-\frac{1}{\rho}\nabla p  + \textbf{g} + \nu \nabla^2 \textbf{u}.\label{eq:ALEmom4}
\end{align}
Replacing operators on the spatial derivative with their $\sigma$-transformed counterparts, now leads to the $\sigma$-transformed equations
\begin{align}
   \frac{\partial \textbf{u}}{\partial \hat{t}}+\frac{\partial \sigma}{\partial t}\frac{\partial \textbf{u}}{\partial \sigma}+ \textbf{u}\cdot\nabla_\sigma\textbf{u} =-\frac{1}{\rho}\nabla_\sigma p  + \textbf{g} + \nu \nabla_\sigma^2 \textbf{u}.\label{eq:ALEsig}
\end{align}
Moreover, (\ref{eq:ALEsig}) is equal to  (\ref{eq:sigmom1}), since for the horizontal momentum equation
\begin{equation}
\begin{aligned}
    \textbf{u}_{\sigma}\cdot  \hat{\nabla}_\sigma u &= u\frac{\partial u}{\partial \hat{x}} + \left(\frac{\partial \sigma}{\partial t} + u\frac{\partial \sigma}{\partial x}  + w\frac{\partial \sigma}{\partial z}\right)\frac{\partial u}{\partial \sigma}\\
    &= \frac{\partial \sigma}{\partial t}\frac{\partial u}{\partial \sigma}+ u\left(\frac{\partial u}{\partial \hat{x}}+\frac{\partial \sigma}{\partial x}\frac{\partial u}{\partial \sigma}\right)+w\left(\frac{\partial \sigma}{\partial z}\frac{\partial u}{\partial \sigma}\right)\\
    &=\frac{\partial \sigma}{\partial t}\frac{\partial u}{\partial \sigma}+\textbf{u}\cdot\nabla_\sigma u,
\end{aligned}
\end{equation}
with the same holding true for the vertical momentum equation. Hence, we conclude that the $\sigma$-transform can be considered a special case of the ALE formulation. 

\section{Temporal discretization}\label{sec:temporal}
To advance the spatially discretized free surface and momentum equations in time, we employ a $s$-stage low-storage explicit Runge-Kutta method (LSERK). Low-storage Runge-Kutta methods are applied to solve initial value problems of the form
\begin{subequations}
\begin{align}
\frac{dy}{dt} &= f(y(t),t), \quad t\geq0, \\
y(t_0)&=y_0, \quad 
\end{align}
\label{eq:ivp}
\end{subequations}
and such LSERK schemes has the general form
\begin{subequations}
\begin{align}
    &y^{(0)} = y^n,\quad K^0 = 0,\\
    &K^k = \alpha_kK^{k-1}+\Delta t f(y^{(k-1)},t_n+c_k\Delta t),\\
    &y^{(k)} = y^{(k-1)}+\beta_kK^k,\quad k=1,2,...,s,\\
    &y^{n+1} = y^{(s)},
\end{align}
\end{subequations}
where $y^n$ refers to the state variables at time $t_n$, $(k)$ refers to the stage number, and $\alpha$, $\beta$ and $c$ are coefficients depending on the chosen LSERK method \cite{Carpenter94}. Note that solely the 5-stage, 4th order LSERK method is used in this work.
\subsection{Conservation of mass}\label{sec:poissonEq}
To ensure that the system has conservation of mass, we have to ensure that $\nabla_\sigma\cdot\textbf{u}^{(k)}=0$ at all stages. This is done using a pressure-correction method where a Poisson boundary value problem is derived for the dynamic pressure that is to be solved at each stage in the Runge-Kutta method to ensure mass conservation through defining a velocity-pressure coupling. We denote the momentum equations in terms of the right hand side function of the initial value problem \eqref{eq:ivp} as
\begin{align}
   f(\textbf{q}) =-\textbf{u}_\sigma\cdot\nabla \textbf{u}-\frac{1}{\rho}(\nabla_\sigma p_D +\hat{\nabla}_\sigma p_S)  + \textbf{g} + \nu \nabla_\sigma^2 \textbf{u},\label{eq:forc1}
\end{align}
with $\textbf{q}=[\textbf{u},p,\eta]$ and $\textbf{u}_\sigma = [u,w_\sigma]$. Combining this with the general LSERK method we get
\begin{align}
    &\textbf{u}^{(k)} = \textbf{u}^{(k-1)}+\beta_{k}K^k
    = \textbf{u}^{(k-1)}+\beta_{k}(\alpha_kK^{k-1}+\Delta t f(\textbf{q}^{(k-1)})).
    \label{eq:rkforc1}
\end{align}
The divergence at stage $k$ can be found by applying the gradient operator (note that $\nabla_\sigma^k$ is the numerical operator at stage $k$, which will be introduced later),
\begin{align}
    \nabla_\sigma^k\cdot\textbf{u}^{(k)} =\nabla_\sigma^k\cdot\textbf{u}^{(k-1)}+\nabla_\sigma\cdot\beta_{k}(\alpha_kK^{k-1}+\Delta t f(\textbf{q}^{(k-1)})).\label{eq:rkdiv1}
\end{align}
Writing out the forcing term and isolating $p_D$ on the LHS, and utilizing that $\nabla^{k}\cdot\textbf{u}^{(k)}=0$, gives the following Poisson equation for the dynamic pressure 
\begin{equation}
\begin{aligned}
    &\nabla_\sigma^k\cdot\nabla_\sigma^{k-1}p_D^{(k-1)} = \frac{\rho}{\beta_k\Delta t}\nabla_\sigma^k\cdot\textbf{u}^{(k-1)}+\frac{\rho\alpha_k}{\Delta t}\nabla_\sigma^k\cdot K^{k-1}\\
    &-\nabla_\sigma^k\cdot\hat{\nabla}_\sigma p_S^{(k-1)}  + \rho\nabla_\sigma^k\cdot \textbf{g} + \rho\nu\nabla_\sigma^k\cdot (\nabla_\sigma^{k-1})^2 \textbf{u}^{(k-1)}-\rho\nabla_\sigma^{k}\cdot(\textbf{u}_\sigma^{(k-1)}\cdot\hat{\nabla}_\sigma\textbf{u}^{(k-1)}),\label{eq:poissonp}
\end{aligned}
\end{equation}
Here we define the {\em mixed-stage Laplace operator} as
\begin{equation}
    \begin{aligned}
        & \nabla^{(k)}_{\sigma} \cdot \nabla^{(k-1)}_{\sigma} = \frac{\partial^2}{\partial x^{*2}}+ \left( \frac{\partial \sigma^{(k-1)}}{\partial x}\frac{\partial \sigma^{(k)}}{\partial x} +\frac{\partial \sigma^{(k-1)}}{\partial z}\frac{\partial \sigma^{(k)}}{\partial z} \right)
        \frac{\partial^2 }{\partial \sigma^2}\\
        & + \left(\frac{\partial \sigma^{(k-1)}}{\partial x} + \frac{\partial \sigma^{(k)}}{\partial x} \right) \frac{\partial^2}{\partial \sigma \partial x^*}  + \frac{\partial^2 \sigma^{(k-1)}}{\partial x^{2}} \frac{\partial }{\partial \sigma}.
    \end{aligned}
    \end{equation}
Remark, \eqref{eq:poissonp} defines a Poisson-type equation for the dynamic pressure $p_D$. To define a boundary value problem for the dynamic pressure, suitable boundary conditions are needed.

\subsection{Boundary conditions for the Poisson BVP}
As mentioned earlier the boundary conditions for the pressure at the surface are given as
\begin{align}
 p_D = 0 \quad \text{on}\quad \Gamma_{FS}.
 \label{eq:bcpd}
\end{align}
For the solid surfaces we utilize the impermeability boundary condition for the velocity to derive a corresponding boundary condition for the pressure. The impermeability boundary condition for the velocity is given as
\begin{align}
\textbf{n}\cdot\textbf{u} = 0\quad \text{on}\quad\Gamma_w\cup \Gamma_b,\label{eq:bcun}    
\end{align}
where $\textbf{n} = (n_x,n_z)^T$ are the normals of the physical domain.\\
To derive the boundary conditions for the Poisson pressure problem, we project (\ref{eq:rkforc1}) in the normal direction to the boundary, and the complete derivation then follows the same ideas as for the conservation of mass in Section \ref{sec:poissonEq}. This results in a bottom boundary condition that need to be fulfilled for consistency in the velocity-pressure coupling defined as 
\begin{equation}
\begin{aligned}
    &\textbf{n}\cdot\nabla_\sigma^{k-1}p_D^{(k-1)} = \frac{\rho}{\beta_k\Delta t}\textbf{n}\cdot\textbf{u}^{(k-1)}+\frac{\rho\alpha_k}{\Delta t}\textbf{n}\cdot K^{k-1}-\textbf{n}\cdot\hat{\nabla}_\sigma p_S^{(k-1)}\\& + \rho\textbf{n}\cdot \textbf{g} + \rho\nu\textbf{n}\cdot (\nabla_\sigma^{k-1})^2 \textbf{u}^{(k-1)}-\rho\textbf{n}\cdot(\textbf{u}_\sigma^{(k-1)}\cdot\hat{\nabla}_\sigma\textbf{u}^{(k-1)}).\label{eq:bcpdneu}
\end{aligned}
\end{equation}
Now, (\ref{eq:poissonp}) in combination with (\ref{eq:bcpd}) and (\ref{eq:bcpdneu}) defines a well-posed Poisson problem that can be solved to determine the dynamic pressure $p_D$ and the divergence free velocity can be updated through the usual LSERK steps. 

\section{Spatial discretization}\label{sec:spatial}
The Navier-Stokes equations along with the resulting Poisson problem for the dynamic pressure are discretized spatially using a nodal spectral element method. The global domain $\Omega$ is partitioned into $N_{el}$ non-overlapping elements $\Omega^n$ such that the union of all elements constitutes the full domain, i.e. $\cup_{n=1}^{N_{el}}\Omega^n=\Omega$. On the domain, we define the finite element space of approximate solutions $V$, consisting of continuous, piece-wise polynomials defined on the discretized elements. More specific, these polynomials are both globally continuous and continuous on each element $\Omega^n$, and are defined as being of degree at most $P$, i.e. $V=\{v\in C^0(\Omega);\forall k \in\{1,...,N_{el}\},v\vert_{\Omega^n}\in \mathbb{P}_P\}$.

\subsection{Kinematic free surface boundary condition}
The weak formulation of the kinematic free surface boundary condition  (\ref{eq:freesurf1}) is found by taking the integral and multiplying with the test function $v(x)$,
\begin{align}
    \int_{\Gamma^n_{FS}}\frac{\partial\eta}{\partial t}v(x) d{\Omega^k_{FS}} = \int_{\Gamma^n_{FS}}\left(-\tilde{u}\frac{\partial\eta}{\partial x} +\tilde{w}\right)v(x)d{{\Gamma^n_{FS}}}.
\end{align}
The equation has to hold for $\eta\in V$ for all $v\in V$. Any function $f(x)$ can be represented globally by piece-wise polynomial functions as
\begin{align}
    f_h = \sum_{i=1}^{K}f_iN_i(x).
\end{align}
Here, $K$ denotes the computational nodes across the entire mesh, while $N_i(x)$ are a set of global finite element basis functions defined such that they have the cardinal property $N_i(x_j)=\delta_{ij}$, with $\delta_{ij}$ being the Kroenecker symbol. This allows us to write the weak formulation as 
\begin{align}
    M\frac{\partial\eta_h}{\partial t} = -A_x^{\tilde{u}_h}+M\tilde{w}_h.
\end{align}
Here we have introduced the global matrices
\begin{subequations}
\begin{align}
    M_{ij} &= \int_{\Gamma^n_{FS}} N_jN_i d{\Gamma^n_{FS}},\\
    (A^b_{x})_{ij} &= \int_{\Gamma^n_{FS}} b(x)\frac{\partial}{\partial x_k}N_jN_i d{\Gamma^n_{FS}}.
\end{align}   %
\end{subequations}

\subsection{Momentum and Poisson equations}\label{sec:laplacewf}
We now define the weak form of the momentum equation (\ref{eq:sigmom1}), meaning we find $\textbf{u}\in V$, such that the following holds for all $v\in V$, 
\begin{align}
   \int_\Omega \frac{\partial \textbf{u}}{\partial t}v d\Omega = \int_\Omega \left(-\textbf{u}_{\sigma}\cdot  \hat{\nabla}_\sigma \textbf{u}- \frac{1}{\rho} (\nabla_{\sigma} p_D + \hat{\nabla}_\sigma p_S) + \textbf{g} + \nu \nabla_{\sigma}^2 \textbf{u}\right)vd\Omega.\label{eq:weak1}
\end{align}
Here integration by parts (IBP) is employed on the Laplacian term due to the second order derivatives. 
For the sake of readability, we define the weak form of the Laplacian by itself,
\begin{equation}
\begin{aligned}
    &\int_\Omega \left(\nabla_{\sigma}^2f\right)v d\Omega =  -\int_{\Omega^k} \frac{\partial v}{\partial x}\frac{\partial f}{\partial x}+\left(v\frac{\partial (\sigma_x^2)}{\partial \sigma}+\frac{\partial v}{\partial \sigma}(\sigma_x^2)\right)\frac{\partial f}{\partial \sigma} 
    \\&+\left(v\frac{\partial (\sigma_z^2)}{\partial \sigma}+\frac{\partial v}{\partial \sigma}(\sigma_z^2)\right)
    \frac{\partial f}{\partial \sigma}+\left(v\frac{\partial \sigma_x}{\partial \sigma}+\frac{\partial v}{\partial \sigma}\sigma_x\right)\frac{\partial f}{\partial x}\\&+\left(v\frac{\partial \sigma_x}{\partial x}+\frac{\partial v}{\partial x}\sigma_x\right)\frac{\partial f}{\partial \sigma}+ v\sigma_{xx}\frac{\partial f}{\partial \sigma} d\Omega^k\\
    &+\int_{\Gamma} v\frac{\partial f}{\partial x}n_x+v(\sigma_x^2)\frac{\partial f}{\partial \sigma}n_\sigma d+v(\sigma_z^2)\frac{\partial f}{\partial \sigma}n_\sigma +v\sigma_x\frac{\partial f}{\partial x}n_\sigma+v\sigma_x\frac{\partial f}{\partial \sigma}n_x d\Gamma.
\end{aligned}\label{eq:weak2}
\end{equation}
The same way we find the weak form of the mixed-stage Poisson problem (\ref{eq:poissonp})
\begin{equation}
\begin{aligned}
    &\int_\Omega\left(\nabla_\sigma^k\cdot\nabla_\sigma^{k-1}p_D^{(k-1)}\right)v d\Omega =\int_\Omega\Big(\Big. \frac{\rho}{\beta_k\Delta t}\nabla_\sigma^k\cdot\textbf{u}^{(k-1)}+\frac{\rho\alpha_k}{\Delta t}\nabla_\sigma^k\cdot K^{k-1}
    \\&-\nabla_\sigma^k\cdot\hat{\nabla}_\sigma p_S^{(k-1)} + \rho\nabla_\sigma^k\cdot \textbf{g} + \rho\nu\nabla_\sigma^k\cdot (\nabla_\sigma^{k-1})^2 \textbf{u}^{(k-1)}-  \\
    &\rho\nabla_\sigma^{k}\cdot(\textbf{u}_\sigma^{(k-1)}\cdot\hat{\nabla}_\sigma\textbf{u}^{(k-1)})\Big.\Big)v d\Omega.
\end{aligned}\label{eq:weak3}
\end{equation}
We apply IBP on the mixed-stage Laplacian,
\begin{equation}
\begin{aligned}
    &\int_\Omega\left(\nabla_\sigma^k\cdot\nabla_\sigma^{k-1}f\right)v d\Omega =  -\int_{\Omega^k} \frac{\partial v}{\partial x}\frac{\partial f}{\partial x}+\left(v\frac{\partial (\overline{\sigma_x^2})}{\partial \sigma}+\frac{\partial v}{\partial \sigma}(\overline{\sigma_x^2})\right)\frac{\partial f}{\partial \sigma}
    \\&+\left(v\frac{\partial (\overline{\sigma_z^2})}{\partial \sigma}+\frac{\partial v}{\partial \sigma}(\overline{\sigma_z^2})\right)
    \frac{\partial f}{\partial \sigma}+\left(v\frac{\partial \sigma_x^{(k-1)}}{\partial \sigma}+\frac{\partial v}{\partial \sigma}\sigma_x^{(k-1)}\right)\frac{\partial f}{\partial x}\\&+\left(v\frac{\partial \sigma_x^{(k)}}{\partial x}+\frac{\partial v}{\partial x}\sigma_x^{(k)}\right)\frac{\partial f}{\partial \sigma}+ v\sigma_{xx}^{(k-1)}\frac{\partial f}{\partial \sigma} d\Omega^k\\
    &+\int_{\Gamma} v\frac{\partial f}{\partial x}n_x+v(\overline{\sigma_x^2})\frac{\partial f}{\partial \sigma}n_\sigma d+v(\overline{\sigma_z^2})\frac{\partial f}{\partial \sigma}n_\sigma +v\sigma_x^{(k-1)}\frac{\partial f}{\partial x}n_\sigma+v\sigma_x^{(k)}\frac{\partial f}{\partial \sigma}n_x d\Gamma.
\end{aligned}\label{eq:weak4}
\end{equation}
We here define special metric coefficients that is a result of the $\sigma$-transform mappings at different stages
\begin{align}
    \overline{\sigma_x^2} =\left(\frac{\partial \sigma^{(k-1)}}{\partial x^*}\frac{\partial \sigma^{(k)}}{\partial x^*}\right),\quad\overline{\sigma_z^2} =\left(\frac{\partial \sigma^{(k-1)}}{\partial z^*}\frac{\partial \sigma^{(k)}}{\partial z^*}\right).
\end{align}
As with in 1D, any function $f(\textbf{x})$ can be represented globally by piece-wise polynomial functions as
\begin{align}
    f_h = \sum_{i=1}^{K}\hat{f}_iN_i(\textbf{x}).\label{eq:fapprox}
\end{align}
Here $K$ denotes the computational nodes across the entire mesh, while $N_i(\textbf{x})$ are a set of global finite element basis functions defined such that they have the cardinal property $N_i(\textbf{x}_j)=\delta_{ij}$, with $\delta_{ij}$ being the Kroenecker symbol. It is possible to represent each of the global basis functions in terms of local basis functions $N^n_i(\textbf{x})$, meaning we can represent $f$ locally on any element as
\begin{align}
    f_h^n=\sum_{j=1}^{N_p}\hat{f}^n_iN^n_i(\textbf{x}),
\end{align}
with $N_p$ being the amount of local nodes on the element. Since $f$ can represent any variable of our governing equations, we only consider the four general cases that occur in our weak formulations (\ref{eq:weak1}-\ref{eq:weak4}),
\begin{align}
\begin{aligned}
    &\int_\Omega vb(\textbf{x})f d\Omega,\quad &\int_\Omega vb(\textbf{x})\frac{\partial f}{\partial x_k} d\Omega,\\
    &\int_\Omega b(\textbf{x})\frac{\partial v}{\partial x_k} \frac{\partial f}{\partial x_k}d\Omega,\quad &\int_\Gamma b(x)v\frac{\partial f}{\partial x_k} nd\Gamma.\label{eq:fweak}
\end{aligned}
\end{align}
We utilize the approximation in (\ref{eq:fapprox}) and choose $v\in\{N_i(\textbf{x})\}_{i=1}^K$ to define a nodal Galerkin scheme. This way the semi-discrete version of (\ref{eq:fweak}) can be written as
\begin{align}
\begin{aligned}
    &M^{b}f,\quad\quad
    &A_{x_k}^{b}f,\\
    &L_{x_k}^b f,\quad\quad
    &B_{x_k}^bf.\label{eq:fsemidiscrete}
\end{aligned}
\end{align}
Here we have introduced the global matrices
\begin{subequations}
\begin{align}
    M^b_{ij} &= \int_\Omega b(\textbf{x})N_jN_i d\Omega,\\
    (A^b_{x_k})_{ij} &= \int_\Omega b(\textbf{x})\frac{\partial}{\partial x_k}N_jN_i d\Omega,\\
    (L^b_{x_k})_{ij} &= \int_\Omega b(\textbf{x})\frac{\partial}{\partial x_k}N_j\frac{\partial}{\partial x_k}N_i d\Omega,\\
    (B^b_{x_k})_{ij} &= \int_\Gamma b(\textbf{x})\frac{\partial}{\partial x_k}N_jN_ind\Gamma.
\end{align}    \label{eq:fglobal}%
\end{subequations}
As we defined the global basis functions with the cardinal property that $N_i(\textbf{x}_j)=\delta_{ij}$, the product of $N_j$ and $N_i$ only contribute if $\textbf{x}_j$ and $\textbf{x}_i$ belong to the same element. We can therefore define local element matrices as 
\begin{subequations}
\begin{align}
    M^n_{ij} &= \int_{\Omega^n} b(\textbf{x})N^n_jN^n_i d\Omega^n,\\
    (A^b_{x_k})^n_{ij} &= \int_{\Omega^n} b(\textbf{x})\frac{\partial}{\partial x_k}N^n_jN^n_i d\Omega^n,\\
    (L^b_{x_k})^n_{ij} &= \int_{\Omega^n} b(\textbf{x})\frac{\partial}{\partial x_k}N^n_j\frac{\partial}{\partial x_k}N^n_i  d\Omega^n,\\
    (B^b_{x_k})^n_{ij} &= \int_{\Gamma^n} b(\textbf{x})\frac{\partial}{\partial x_k}N^n_jN^n_ind\Gamma^n.
\end{align}    \label{eq:flocal}%
\end{subequations}
The local element matrices are essentially the local contribution to the global matrices in (\ref{eq:fglobal}). Due the the domain partition introduced earlier, we can therefore define the global sparse matrices as a sum of the local dense element matrices
\begin{align}
    M_{ij}=\sum_{n=1}^{N_{el}}M^n_{ij}=\sum_{n=1}^{N_{el}}\int_{\Omega^n} N^n_jN^n_i d\Omega^n.\label{eq:globallocal}
\end{align}
The other global matrices can be constructed similarly. 

\subsection{Element construction in 2D}
To represent the solution in the 2D space we consider quadrilateral elements. We introduce a reference element given by $\mathcal{T}=\{(r,t)\in \mathbb{R}:-1<(r,t)<1\}$, on which we will define element basis functions and node positions. It turns out defining these on the reference element and then mapping to the general elements is an effective way of handling the representation. We construct the 2D basis as a tensor product of the 1D orthonomal Jacobi polynomials given by $\tilde{P}^{(\alpha,\beta)}_k(x)$ on the interval $x \in [-1,1]$. Here $k$ denotes the arbitrary order, and $\alpha$ and $\beta$ are the parameters. Through this we can define the element 2D basis as
\begin{align}
    \psi_{nm}(r,s) = \tilde{P}^{(0,0)}_n(r)\tilde{P}^{(0,0)}_m(s),  
\end{align}
with $\alpha=\beta=0$ leading to the Legendre polynomials, which are a special case of the Jacobi polynomials and will be denoted as $\tilde{P}_k$ from now on. The Legendre polynomials can efficiently be computed through a recurrence relation given as
\begin{subequations}
\begin{align}
    a_{k}\tilde{P}_{k}(x) &= x\tilde{P}_{k-1}(x) - a_{k-1}\tilde{P}_{k-2}(x),\\
    a_k &= \sqrt{\frac{k^2}{(2k+1)(2k-1)}}, 
\end{align}
\end{subequations}
where the first two Legendre polynomials are defined as
\begin{align}
    \tilde{P}_0(x) = \frac{1}{\sqrt{2}}, \qquad \tilde{P}_1(x) = \sqrt{\frac{3}{2}}x.
\end{align}
The nodal distribution on each element is given by the Legendre-Gauss-Lobatto points. The modal basis functions defined through Legendre polynomials all have corresponding Lagrange polynomials, which means a continuous function $f_h$ can be represented locally in the form of either a modal or nodal expansion
\begin{align}
    f_h(\textbf{r}) \approx \sum_{m=1}^{N_{p}}\hat{f}_m^n\psi(\textbf{r})=\sum_{m=1}^{N_{p}}f_m^nh_m(\textbf{r}).
\end{align}
The generalized Vandermonde matrix $\mathcal{V}$ and its derivative variants can be defined as
\begin{align}
\mathcal{V}_{ij} = \psi_j(\textbf{r}_i),\quad
    (\mathcal{V}_r)_{ij} = \partial_r\psi_j(\textbf{r}_i),\quad (\mathcal{V}_s)_{ij} = \partial_s\psi_j(\textbf{r}_i).
\end{align}
From this, the relation between the modal and nodal expansion can be written as
\begin{align}
    \textbf{f}_h = \mathcal{V}\hat{\textbf{f}},
\end{align}
and we can construct the Lagrange polynomials with the cardinal property $h_i(x_j)=\delta_{ij}$, along with their derivatives, as 
\begin{subequations}
\begin{align}
    h_i(\textbf{r})&=\sum_{j=1}^{N_p}(\mathcal{V}^T)^{-1}_{ij}\phi_j(\textbf{r}),\\
    \partial_r h_i(\textbf{r})&=\sum_{j=1}^{N_p}(\mathcal{V}^T)^{-1}_{ij}\partial_r\phi_j(\textbf{r}) = \sum_{j=1}^{N_p}(\mathcal{V}^T)^{-1}_{ij}(\mathcal{V}_r)_{ij}\phi_j(\textbf{r}),\\
    \partial_s h_i(\textbf{r})&=\sum_{j=1}^{N_p}(\mathcal{V}^T)^{-1}_{ij}\partial_s\phi_j(\textbf{r}) = \sum_{j=1}^{N_p}(\mathcal{V}^T)^{-1}_{ij}(\mathcal{V}_s)_{ij}\phi_j(\textbf{r}).
\end{align}    
\end{subequations}
From these we can construct the local mass matrix, along with local derivative operators as
\begin{align}
\begin{aligned}
    \mathcal{M}_{ij} = (\mathcal{V}\mathcal{V}^T)^{-1}_{ij},\quad
    \mathcal{D}_r = \mathcal{V}_r\mathcal{V}^{-1},\quad
    \mathcal{D}_s = \mathcal{V}_s\mathcal{V}^{-1}.
\end{aligned}
\end{align}
To obtain the derivative operators in physical space, we let $r_x = \partial_x r$ and $s_z = \partial_z s$, and define $R_x = \text{diag}(r_x)$ and $S_z = \text{diag}(s_z)$. We can then apply the chain rule, resulting in the physical derivative operators 
\begin{align}
\begin{aligned}
    &\mathcal{D}_x = R_x\mathcal{D}_r,\quad
    &\mathcal{D}_z = S_z\mathcal{D}_s.
    \end{aligned}
\end{align}
With these operators, the local element matrices defined in (\ref{eq:flocal}) can be computed from the reference element operators as
\begin{subequations}
\begin{align}
    M^n & = \lvert\mathcal{J}^n\rvert \mathcal{M},\\
    (A^b_{x_k})^n &=\lvert\mathcal{J}^n\rvert D_{x_k}M\mathcal{I},\\
    (L^b_{x_k})^n &= \lvert\mathcal{J}^n\rvert D_{x_k}MD_{x_k},\\
    (B^b)^n &= \sum_i \lvert(\mathcal{J}^{s})^n_i\rvert D^{1D}_{x^k} M^{1D}\mathcal{I}.
\end{align}    \label{eq:flocal2}%
\end{subequations}
Here $\mathcal{J}^n$ is the Jacobian of the affine mapping $T^n:\Omega^n\rightarrow\Omega^r$, where $\Omega^r$ is the reference element, $(\mathcal{J}^s)^n_i$ is the corresponding surface Jacobian for the $i$'th surface, and $\mathcal{I}$ is the identity matrix. 

\subsection{Spectral Filtering}\label{sec:secDamp}
The strong nonlinear terms present in the governing equations can present a challenge for maintaining stability in the simulation, as nonlinear waves can develop unstable sawtooth waves when propagated over time. To avoid this issue, we employ a spectral filtering strategy, taking advantage of the dual nodal-modal representation of the variables. We define a exponential cut-off filter as 
\begin{equation}
     S(i) = \begin{cases} 
          1 & 0\leq i\leq P_c, \\
          \exp\left(\alpha\left(\frac{i-P_c}{P+1-P_c}\right)^\beta\right) & P_c<i\leq P. 
       \end{cases}
\end{equation}
The filter is applied to both the surface elevation $\eta$ and the velocities $u,w$, with the filtered local element solution of any given variable in modal representation defined as
\begin{align}
    f^n(\textbf{x})=\sum_{m=1}^{N_{p}}S(i)\hat{f}_m^n\psi(T^n(\textbf{x})).
\end{align}
The filter reduces energy in any modes above the cutoff $P_c$, with parameters $\alpha$ and $\beta$ determining the strength of the filtering. To maintain the spectral properties of the model, the parameters are generally chosen such that only a few percent of the energy of the highest modes are removed.
The modes can easily be obtained through use of the Vandermonde matrices, with the full filtering operation gives as
\begin{align}
    \mathcal{F} = \mathcal{V}F\mathcal{V}^{-1},\quad F=\text{diag}(S).
\end{align}

\subsection{Transforming the outward pointing normals}
Due to how the surface integrals appear in the weak form of the $\sigma$-transformed operators, the boundary conditions need to be defined in terms of the outward pointing normals in the $\sigma$-domain. The outward pointing normals can be transformed as
\begin{align}
    \begin{pmatrix}
n_x\\
n_z
\end{pmatrix}=\left(J^{T}\begin{pmatrix}
n_x^*\\
n_\sigma
\end{pmatrix}\right)*\|J^{T}\begin{pmatrix}
n_x^*\\
n_\sigma
\end{pmatrix}\|_2^{-1}.\label{eq:outVecTrans1}
\end{align}
If we let
\begin{align}
    N = \|J^{T}\begin{pmatrix}
n_x^*\\
n_\sigma
\end{pmatrix}\|_2,
\end{align}
we can rewrite (\ref{eq:outVecTrans1}) as
\begin{align}
    N\begin{pmatrix}
n_x\\
n_z
\end{pmatrix}=\left(J^{T}\begin{pmatrix}
n_x^*\\
n_\sigma
\end{pmatrix}\right).
\end{align}
We combine this with the LHS of our boundary condition
\begin{equation}
\begin{aligned}
    &Nn\cdot\nabla_\sigma p = \left(J^{T}\begin{pmatrix}
n_x^*\\
n_\sigma
\end{pmatrix}\right)\cdot\nabla_\sigma p= \left(\begin{bmatrix}
1 & \frac{d\sigma}{dx}\\
0 & \frac{d\sigma}{dz}
\end{bmatrix}\begin{pmatrix}
n_x^*\\
n_\sigma
\end{pmatrix}\right)\cdot\begin{pmatrix}
\frac{dp}{dx^*}+\frac{d\sigma}{dx}\frac{dp}{d\sigma}\\
\frac{d\sigma}{dz}\frac{dp}{d\sigma}
\end{pmatrix}\\
&=n_x^*\frac{dp}{dx^*}+n_\sigma\left(\frac{d\sigma}{dx}\right)^2\frac{dp}{d\sigma}+n_\sigma\left(\frac{d\sigma}{dz}\right)^2\frac{dp}{d\sigma}+n_\sigma\frac{d\sigma}{dx}\frac{dp}{dx^*}+n_x^*\frac{d\sigma}{dx}\frac{dp}{d\sigma}\label{eq:surfMatch1}.
\end{aligned}
\end{equation}
The terms on the RHS matches the terms in the surface integrals of the weak form in section \ref{sec:laplacewf}, meaning that by rewriting the boundary conditions to the form
\begin{align}
    N*(n\cdot\nabla_\sigma p) = N*(g(x,\sigma,t)),
\end{align}
we can enforce boundary conditions through the boundary terms in the weak formulation in the usual manner.

\section{Iterative geometric $p$-multigrid method}\label{sec:mg}
\begin{figure}[H]
    \centering
    \includegraphics[width=0.6\linewidth]{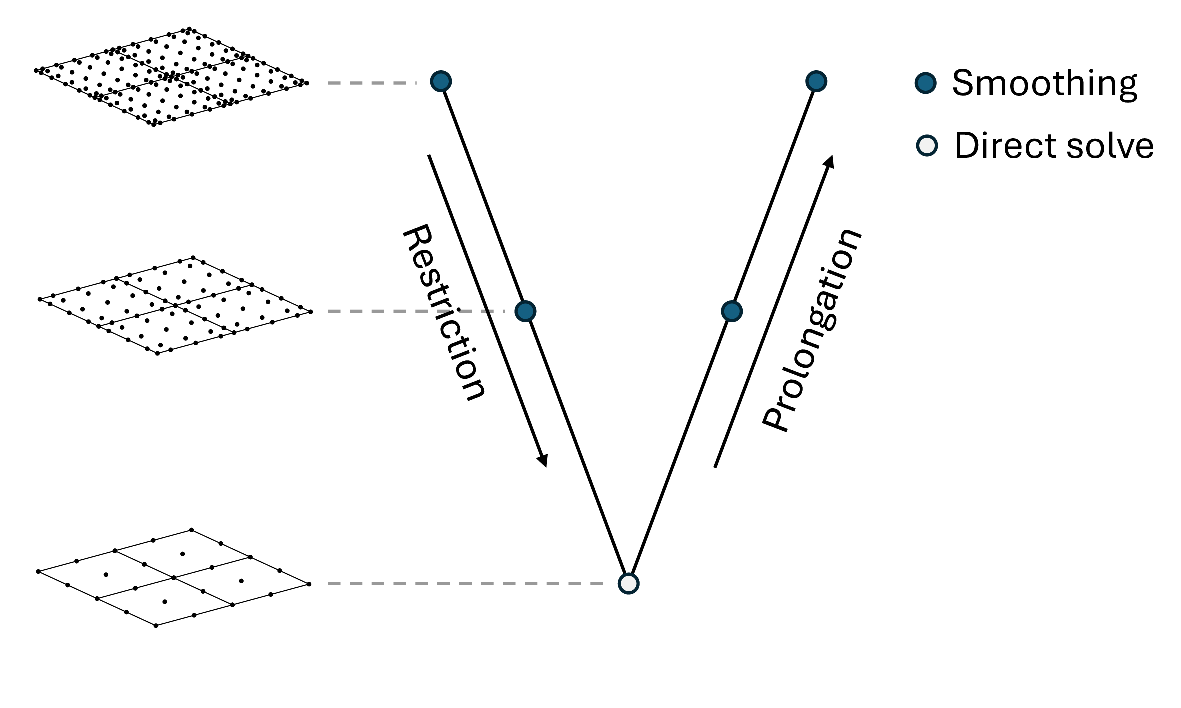}
    \caption{Example of $V$-cycle for geometric $p$-multigrid on a $2\times2$ mesh.}
    \label{fig:vcycle}
\end{figure}

The Poisson problem for the dynamic pressure has to be solved at every stage of the RK method, and an efficient way to do this is required. The class of multigrid methods \cite{Trottenberg01} takes advantage of the performance of basic stationary iterative solvers, which are great at eliminating high-frequency errors but less so with lower-frequency ones. By transferring errors to a coarser grid, lower-frequency errors are converted into high-frequency ones, allowing the basic iterative method to be effective once more. Using a high-order spectral method, such as SEM, allows the use of the geometric $p$-multigrid method, which takes advantage of the fact that the coarsening of the grid can be done by simply reducing the polynomial order while keeping the mesh topology unchanged. In this paper, we will restrict the focus to the $V$-cycle, which is illustrated in Figure \ref{fig:vcycle}. Starting at the finest grid, the error of the solution is transferred to coarser and coarser grids, applying smoothing at each of these, until reaching the coarsest grid. Then, the estimated error is transferred back up through the grid levels to the finest grid, again applying smoothing at each level.

\subsection{Grid strategy}
An integral part of the geometric $p$-multigrid method is establishing the grid hierarchy. As mentioned, this is done by reducing the polynomial order $p$ of the SEM discretization. We let $G^N$ be the original fine grid of order $P_N$ and $G^1$ be the coarsest grid of order $P_1$. The subsequent coarser grids are then generated by reducing the order as
\begin{align}
    P_{n-1} = \text{ceil}\left(\frac{P_{n}+1}{2}\right). \label{eq:gridcoars}
\end{align}
For cases where the horizontal and vertical order of discretization differs, a semi-coarsening strategy is employed, reducing the highest order of discretization until both directions have matching orders.  

\subsection{Transfer operators}
A necessity for the geometric $p$-multigrid method is an efficient way to transfer the error between grids. This is done through the use of restriction operators moving the errors from fine to coarse, denoted as $R_{fc}$, and prolongation operators moving the errors from coarse to fine, denoted as $P_{cf}$. In the case of SEM, the modal representation lends itself well to creating linear operators through interpolation. Locally, the solution can be represented on grid $P$ as
\begin{align}
    f_P^n = \mathcal{V}\hat{f}_P^n,
\end{align}
where it follows that the solution can be represented on grid $p$ as
\begin{align}
    f_p^n = \Phi^p\hat{f}_P^n,
\end{align}
with $\mathcal{V}_{ij} = \psi_j(\textbf{x}^P_i)$ and $\Phi^p_{ij} = \psi_j(\textbf{x}^p_i)$. From these two representations, it follows that any solution can be transferred from grid $P$ to grid $p$ through
\begin{align}
    f_p^n = \Phi^p\mathcal{V}^{-1} f_P^n,
\end{align}
meaning we can define the transfer operator as
\begin{align}
    \mathcal{I}_{Pp} = \Phi^p\mathcal{V}^{-1}.
\end{align}
This means that the prolongation operator is simply defined as
\begin{align}
    P_{cf} = \mathcal{I}_{fc}.
\end{align}
While the restriction operator could be defined in the same manner, we instead define it as the transpose of the prolongation operator \cite{Rønquist1987},
\begin{align}
    R_{fc} = P_{cf}^T.
\end{align}

\subsection{Smoothing operator}
The smoothing operator consists of a stationary iterative method which serves to reduce the error at every grid level. An important note is that the iterative method is not applied until convergence, but rather just applied for a given amount of iterations. As mentioned earlier, this allows for fast reductions of high-frequency errors at every grid level. In this work, we apply the additive Schwarz method (ASM) as the smoother. This method has been shown to be an effective smoother for higher order methods such as SEM \cite{AllanLaskowski,Lottes05}. We can represent the smoothing process as 
\begin{align}
    \textbf{x}^{k+1} = \textbf{x}^k-\mathcal{S}^{-1}(\mathcal{A}\textbf{x}^k-\textbf{b}),
\end{align}
with $\mathcal{A}$ being the system matrix, and $\mathcal{S}^{-1}$ is the Schwarz preconditioner. To construct the preconditioner we define transfer operators $\mathcal{R}_i: \Omega_k \rightarrow \hat{\Omega}_k$ from the original computational domain, to a new subdomain $\hat{\Omega}_k$ defined as an element-block combined with a user-determined amount of overlap nodes with neighbouring blocks. This allows us to construct the Schwarz preconditioner as an element-wise block matrix given as
\begin{align}
    \mathcal{S}^{-1}=\mathcal{W}\sum_{i=1}^{N_{el}}\mathcal{R}_i^T\left(\mathcal{R}_i\mathcal{A}\mathcal{R}_i^T\right)^{-1}\mathcal{R}_i.
\end{align}
Here $\mathcal{W}$ is a weight operator defined as
\begin{align}
    \mathcal{W}=\left(\sum_{i=1}^{N_{el}}\mathcal{R}_i\mathcal{R}_i^T\right)^{-1}.
\end{align}

\subsection{Combination with stationary iterative solver}
On non-moving meshes, the cost of computing the geometric multigrid operators is usually negligible, as this setup stage is only required at the start of any given simulation. However, due to the movement of our mesh, represented through the $\sigma$-transform, these operators would in principle have to be recomputed at every RK stage. To circumvent this issue, we combine the multigrid solver with both a preconditioned defect-correct method (PDC) \cite{EngsigKarup2014} and GMRES \cite{SS86}, utilizing the multigrid solver to solve the occurring pre-conditioning problem
\begin{align}
    \mathcal{M}^{-1}\textbf{e} = -\textbf{r}.
\end{align}
By basing the pre-conditioning matrix $\mathcal{M}$ on a linearised Navier-Stokes formulation arising from introducing the small--amplitude wave assumption $H/L\ll \mathcal{O}(1)$ to construct a time-constant preconditioner, it is only necessary to compute the multigrid operators at the beginning of any given simulation, thus making the computational cost of the setup stage negligible. 

\section{Numerical experiments}\label{sec:res}
\subsection{Model verification using nonlinear streamfunction waves}
\begin{figure}[H]
\vspace*{0cm}
\hspace*{0cm}
\centering
\subfloat[]{
  \includegraphics[width=38mm]{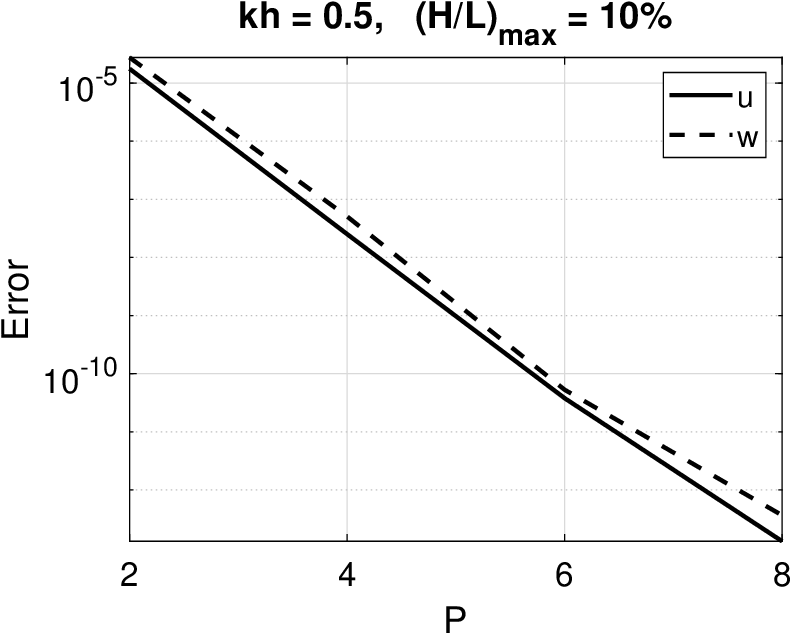}
}
\subfloat[]{
  \includegraphics[width=38mm]{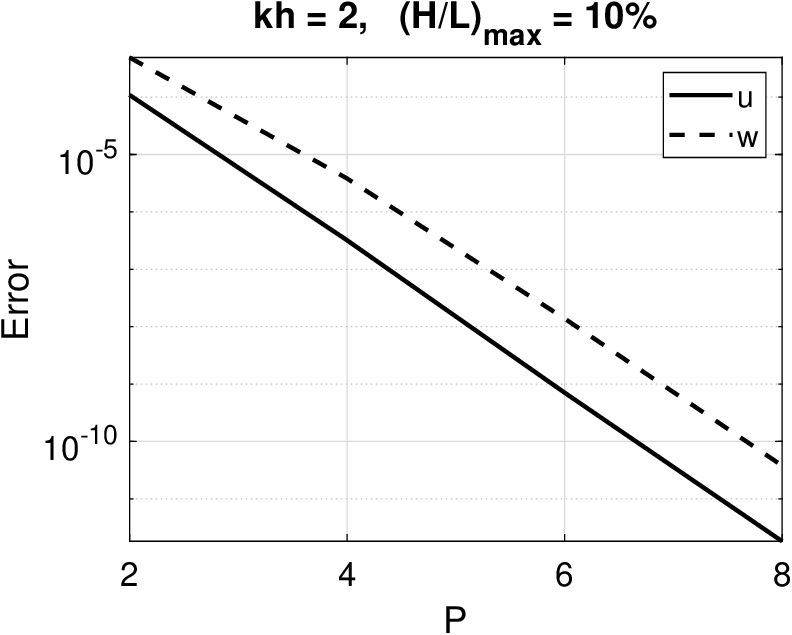}
}
\subfloat[]{
  \includegraphics[width=38mm]{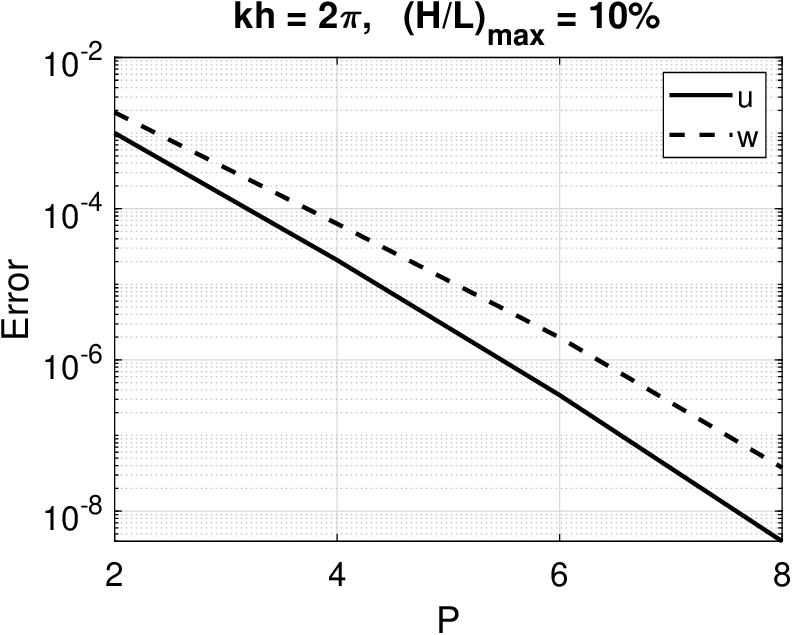}
}
\hspace{0mm}
\hspace*{0cm}
\subfloat[]{
  \includegraphics[width=38mm]{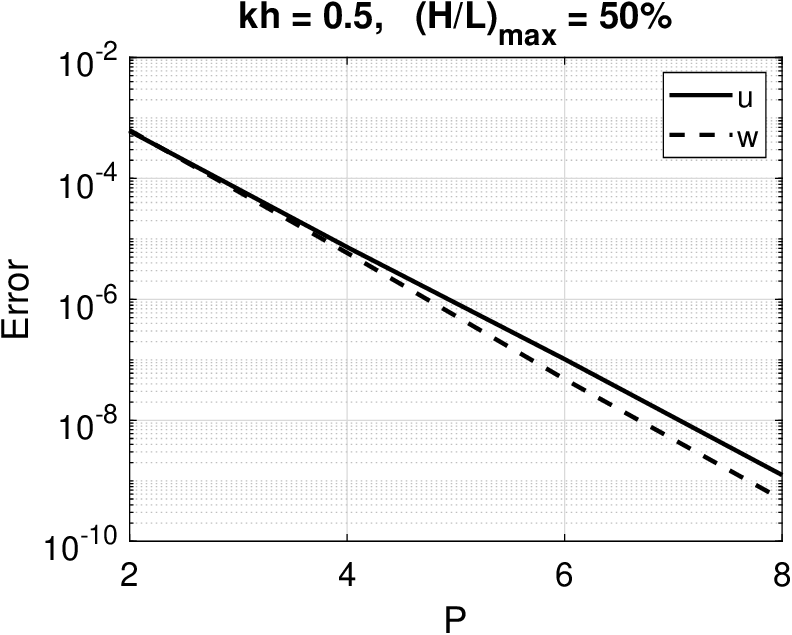}
}
\subfloat[]{
  \includegraphics[width=38mm]{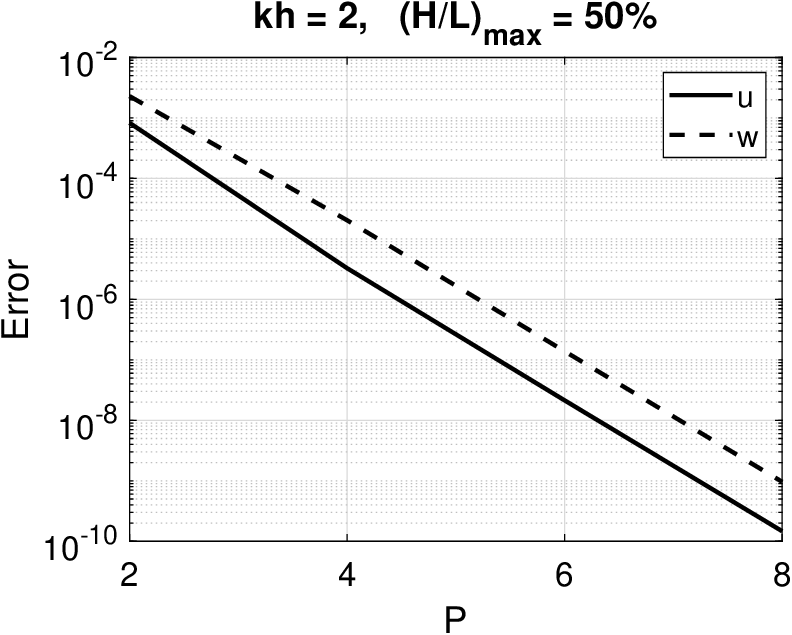}
}
\subfloat[]{
  \includegraphics[width=38mm]{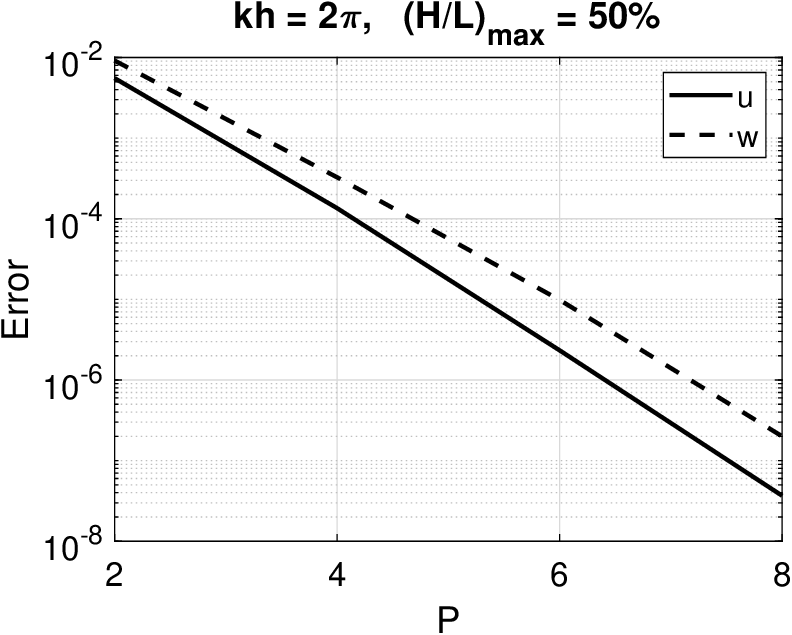}
}
\hspace{0mm}
\hspace*{0cm}
\subfloat[]{
  \includegraphics[width=38mm]{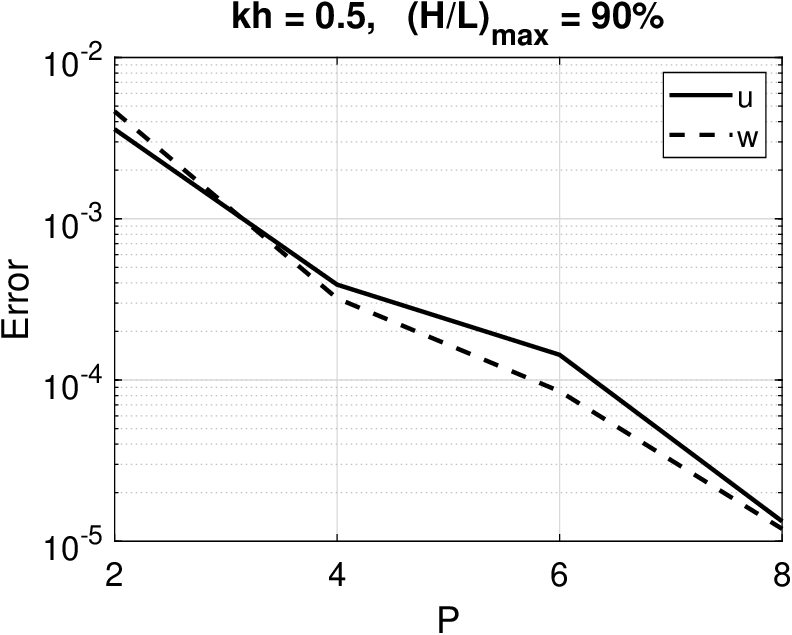}
}
\subfloat[]{
  \includegraphics[width=38mm]{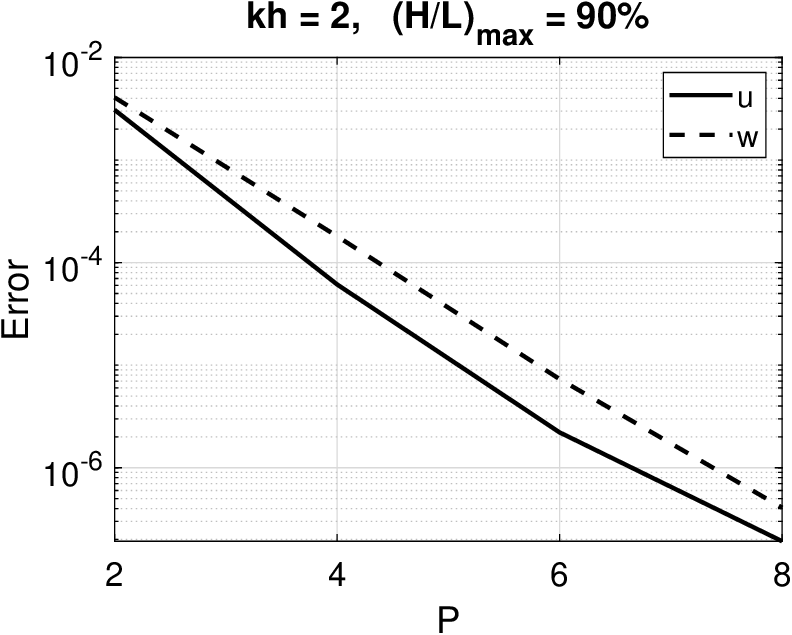}
}
\subfloat[]{
  \includegraphics[width=38mm]{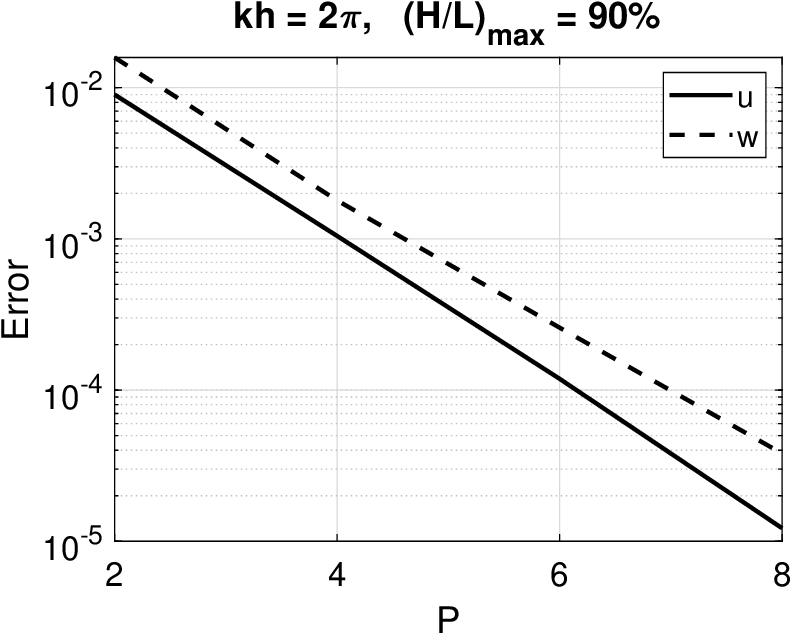}
}
\caption{Error convergence for the velocities $u$ and $w$ across varying depth and wave steepness.}
\label{fig:convvel}
\end{figure}

In order to verify that the model is capable of achieving the expected spectral $p$-convergence, we here present a convergence study of the error in approximation of the velocity variables. We here assume no viscosity i.e. we let $\nu = 0\text{ }[\frac{m^2}{s}]$, which allows us to compare results against analytical stream function solutions \cite{Rienecker81} which are valid solutions for non-viscous and irrotational fluid flow. We use a fixed number of elements in both the vertical and horizontal direction, and only vary the polynomial order $P$ to demonstrate $p$-convergence. To allow the study to take into account the error from both the Poisson pressure problem, the momentum equations and the free surface equation, the error is computed after one full time step of the RK method, and is defined as
\begin{align}
    \text{Error} = \|u-u_e\|_\infty,
\end{align}
with $u_e$ being the exact streamfunction solution. The amount of elements are kept constant at $Ne_x = 20, Ne_y=2$, and the time step size is chosen such that spatial errors dominate. Figure \ref{fig:convvel} shows the error convergence for nine different cases. These are shallow ($kh=0.5$), intermediate ($kh=2$) and deep ($kh=2\pi$) water with three cases of increasingly steep waves for each depth. Note that we define wave steepness as a percentage of the maximum allowed steepness before breaking, as defined by Battjes \cite{Battjes1},
\begin{align}
    \left(\frac{H}{L}\right)_{max} = 0.1401 \tanh(0.8863kh).
\end{align}

\subsection{Harmonic generation over a submerged bar}\label{sec:bar}
\begin{figure}[H]
    \centering
    \includegraphics[width=1\linewidth]{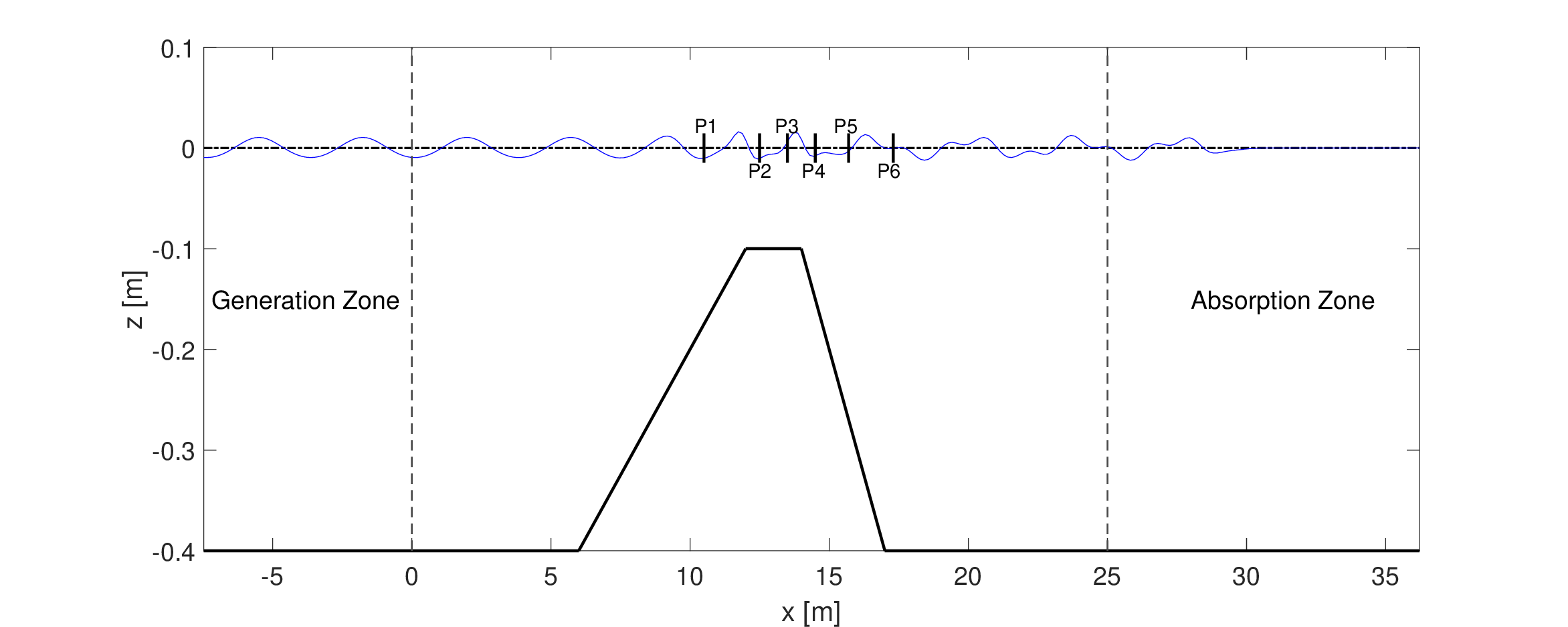}
    \caption{Setup of the bar test.}
    \label{fig:bartest}
\end{figure}

To show that the solver is capable of handling both wave-tank setups, as well as non-flat bathymetry, we present results of harmonic generation over a submerged bar, which was performed both experimentally and numerically by Beji and Battjes \cite{Beji94}. The setup of the experiment can be seen in figure \ref{fig:bartest}, and consists of a raised bar in the middle of a wave tank. Waves are generated in the generation zone and propagated across the bar, before being absorbed in the absorption zone. The zones are defined as proposed by Larsen and Dancy \cite{Larsen83}, with the full definitions available in Appendix \ref{app:waveZones}.
The slope of that bar incurs a shoaling effect which steepens the incoming waves, before they start decomposing into free wave harmonics after the top. This results in a rapidly changing wave profile due to  wave dispersion, that is impossible to capture if the model does not accurately describe the wave phenomena. 

Mildly nonlinear waves are initially generated with height of $H=0.02 \text{m}$, length of $L = 3.74$ m and period of $T = 0.02$ s giving a dimensionless depth of $kh = 0.67$. The domain is with elements of order $P = 8$, which are evenly spaced such that $N_x = 100$ and $N_z = 2$ for the horizontal and vertical directions respectively.

\begin{figure}[H]
\hspace*{0cm}
\centering
\subfloat[]{
  \includegraphics[width=60mm]{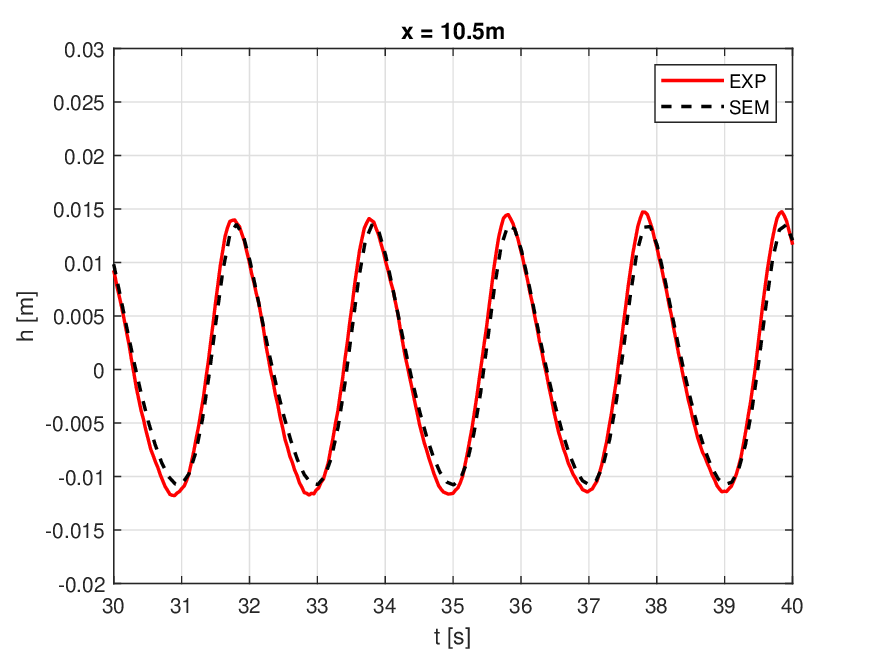}
}
\subfloat[]{
  \includegraphics[width=60mm]{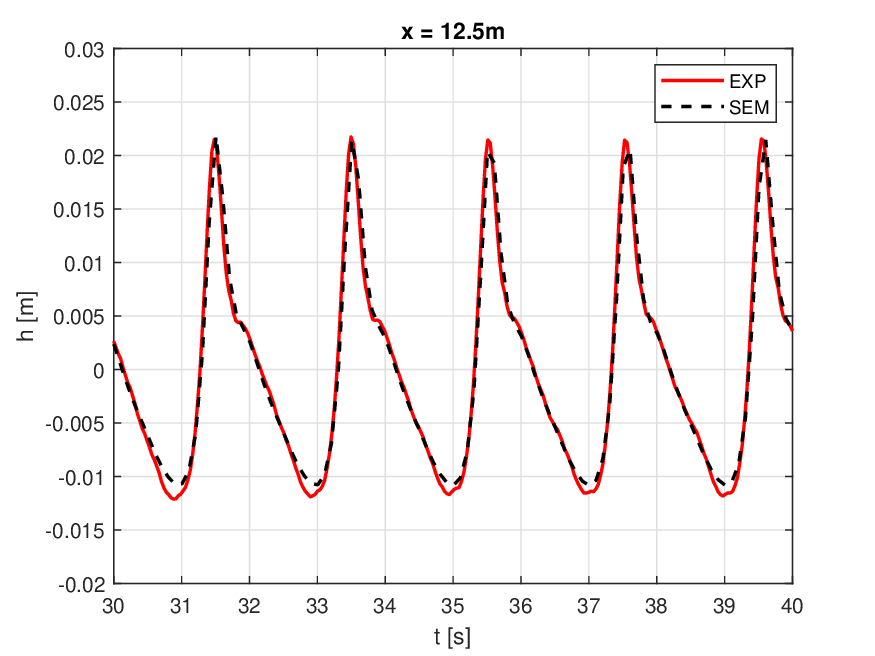}
}

\hspace*{0cm}
\subfloat[]{
 \includegraphics[width=60mm]{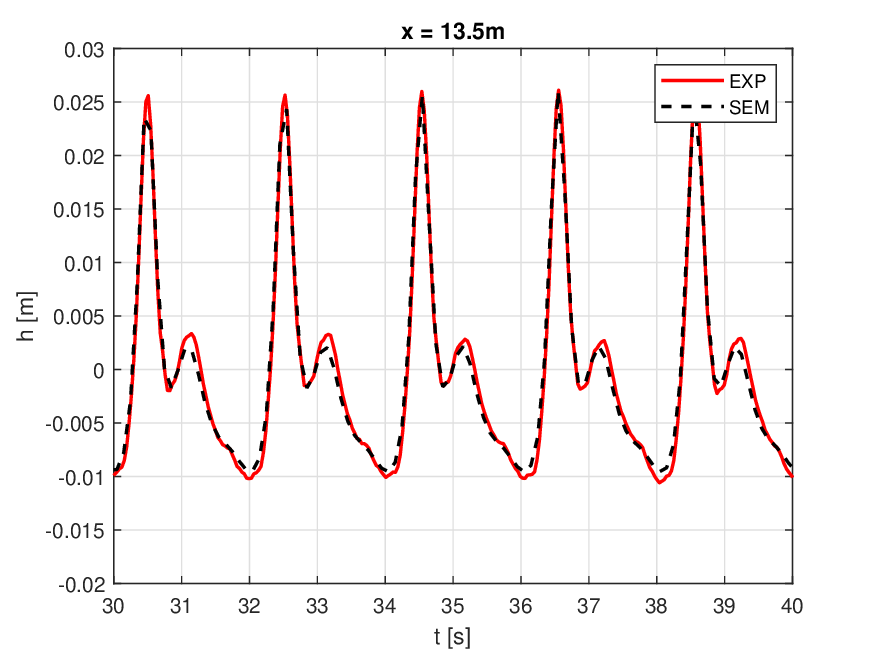}
}
\subfloat[]{
  \includegraphics[width=60mm]{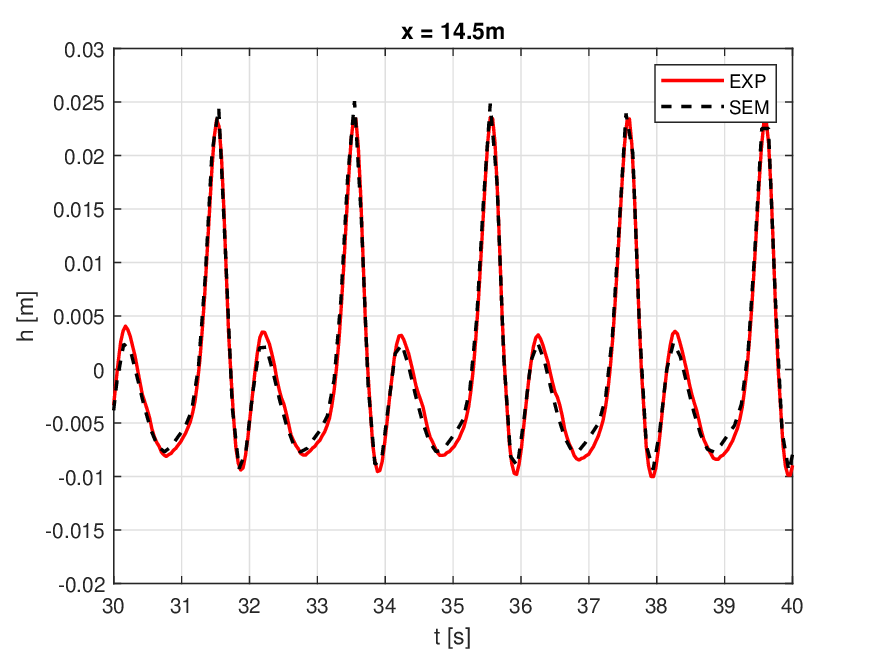}
}

\hspace*{0cm}
\subfloat[]{
  \includegraphics[width=60mm]{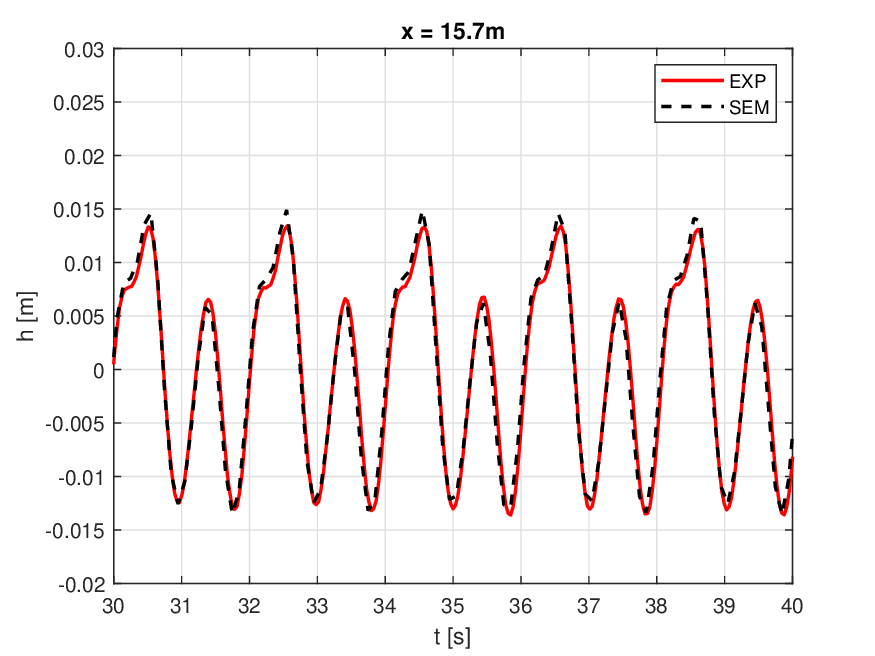}
}
\subfloat[]{
  \includegraphics[width=60mm]{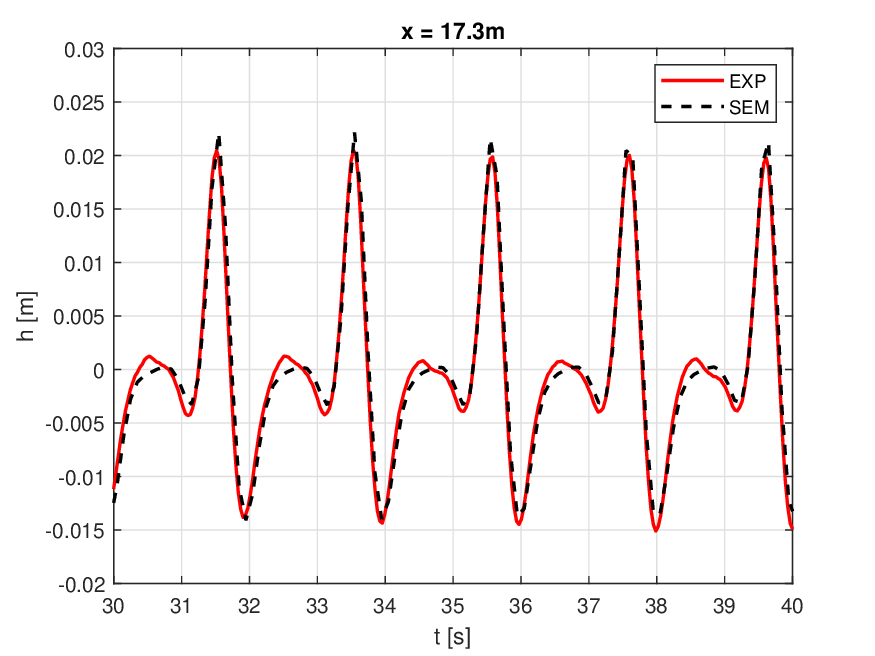}
}
\caption{Comparison of numerical and experimental results for a harmonic generation over a submerged bar.}
\label{fig:subBar2}
\end{figure}

From figure \ref{fig:subBar2}, it is clear that the model is capable of capturing the nonlinear effects as seen by the excellent agreement between simulated and experimental results. 

\subsection{Geometric $p$-multigrid efficiency}
\begin{figure}[H]
\hspace*{-0.8cm}
\centering
\subfloat[]{
  \includegraphics[width=65mm]{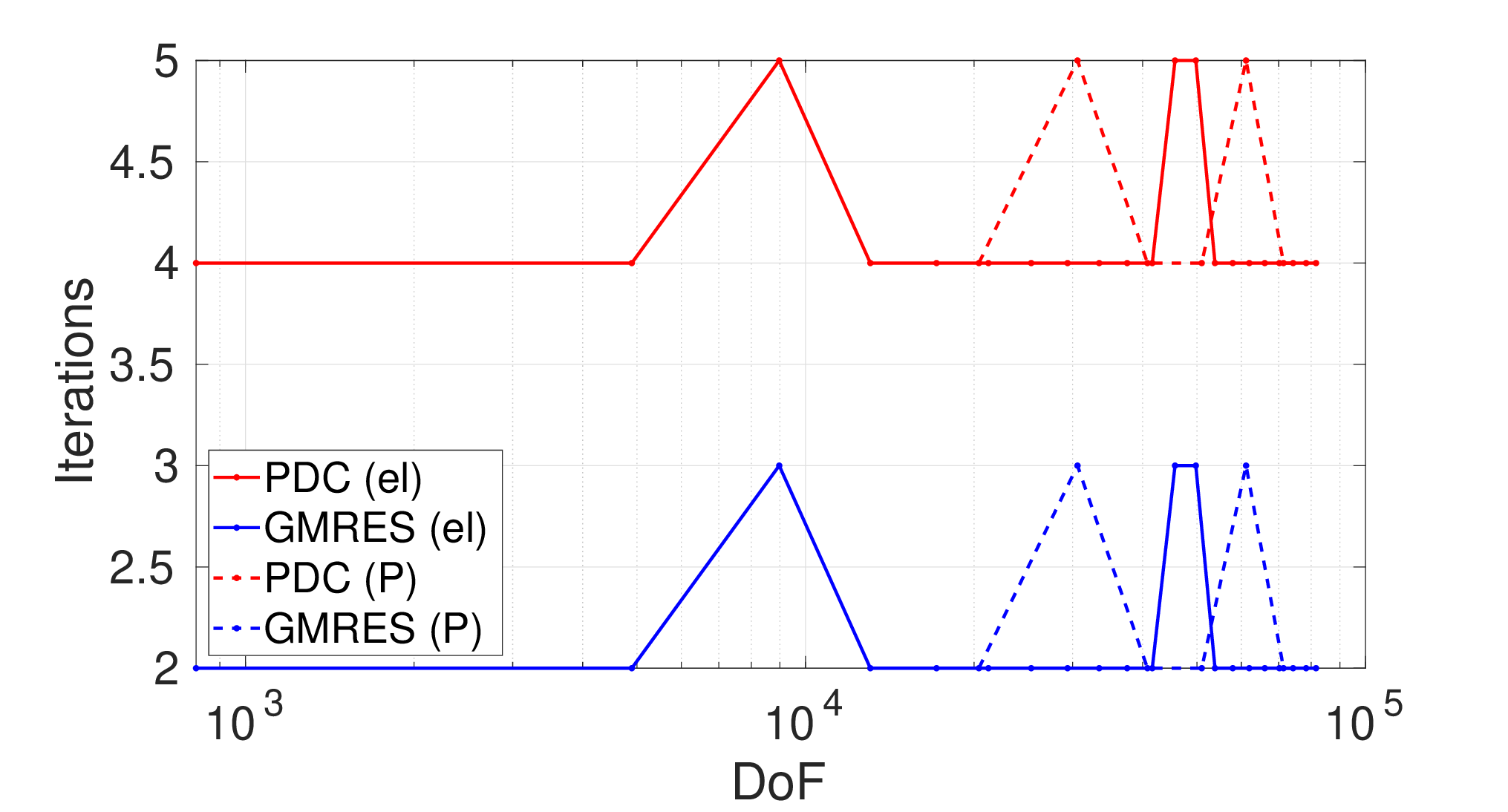}
}
\subfloat[]{
  \includegraphics[width=65mm]{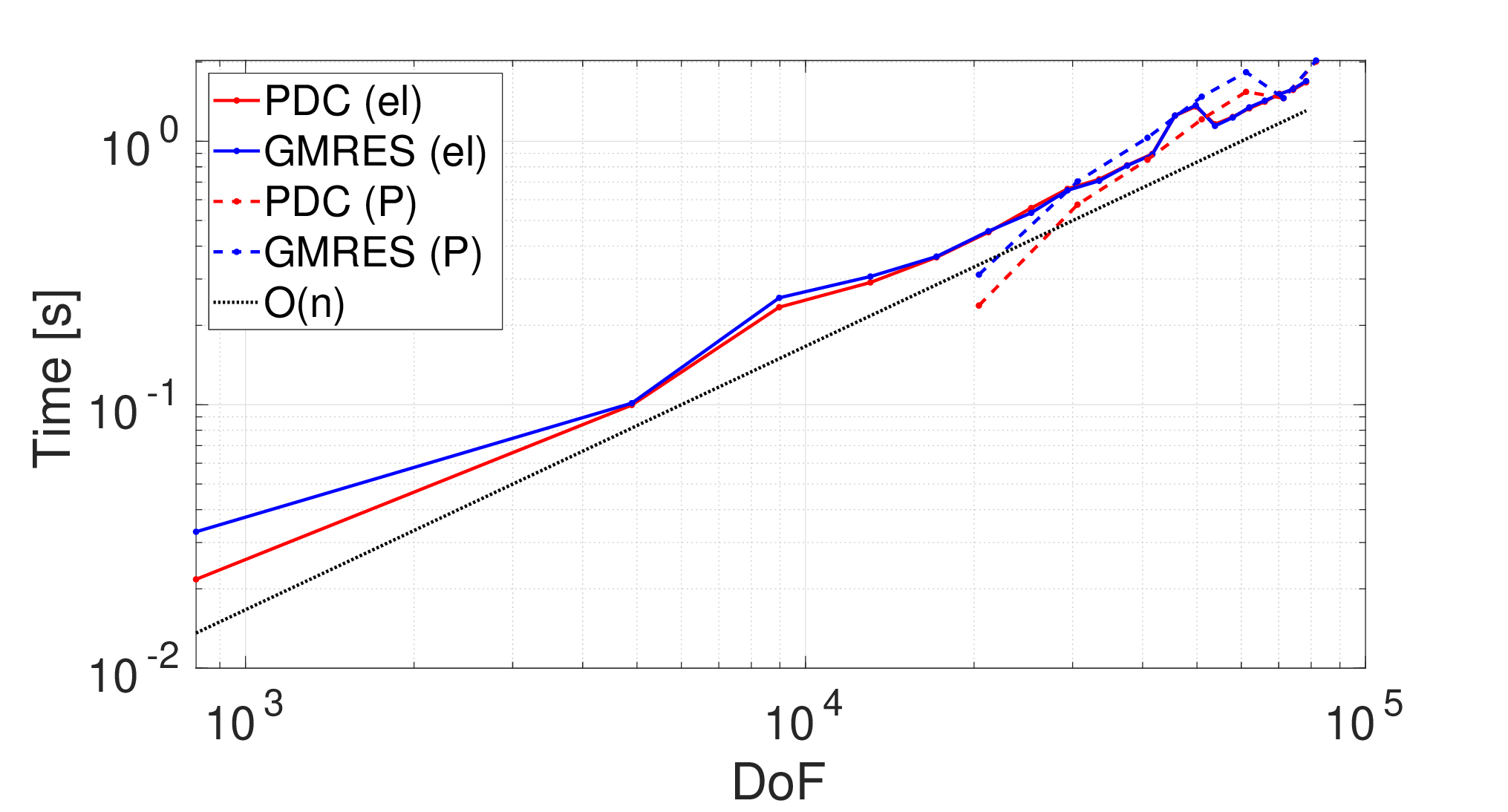}
}
\caption{Computational scaling of $p$-multigrid accelerated PDC and GMRES.}
\label{fig:mgscaling}
\end{figure}

We here present results to show the efficiency of the geometric $p$-multigrid accelerated PDC method. Table \ref{tab:ittimemg} shows iterations and computational cost of solving the occurring Poisson pressure problem for a discretization at roughly the same degrees of freedom as the bar test with discretization parameters $P = 8$, $N_x = 102$ and $N_z = 2$. The simulation is undertaken for seventeen waves of $kh=1$, and wave steepness of $(H/L)_{\text{max}}=30\%$ corresponding to $(H/L)=0.0301$, leading to $\sim 48$ points per wavelength (PPW). The computational cost is compared to the built-in direct solver in Matlab R2021a.
\begin{table}[h]
\centering
\begin{tabular}{|l|l|l|l|}
\hline
       Method             & $10^{-4}$ & $10^{-6}$ & $10^{-8}$ \\ \hline
PDC-MG Iterations   &          3                     &               5                 &                8                \\ 
PDC-MG time [s] &            0.1402                    &                0.2736                &                  0.4079              \\ \hline
GMRES-MG Iterations &            1                    &                3                &                  4              \\ 
GMRES-MG time [s] &            0.1614                    &                0.2721                &                  0.3589              \\\hline 
Direct time [s] &                   0.1643         & 0.1643& 0.1643                                                                  \\ \hline
\end{tabular}
\caption{Iterations and computational cost to be within specified error tolerance.}\label{tab:ittimemg}
\end{table}
While the solvers slightly outperforms the direct solver for high tolerances, the direct solver tends to outperform both the PDC-MG and GMRES-MG methods at lower tolerances. However, keep in mind that iterative solvers tend to be outperformed by direct solvers at lower degrees of freedom (DoF) and when the sparse system matrix has a  small band-structure as is the case for a setup with few elements in the vertical  and in two space dimensions. It is therefore of interest to look at the computational scaling properties of the solvers. We consider the scaling both in terms of increasing DoF through the number of horizontal elements $N_x$ and the horizontal order $P_x$. The vertical discretization is fixed at $N_z=2$ and $P_z=8$, and the Poisson problem is solved for waves of $kh=1$ and $(H/L)_{\text{max}}=30\%$ with a fixed PPW at $\sim 48$. For the tests increasing the number of elements $N_x$, we fix the order at $P_x=8$, and for the tests increasing the order of discretization $P_x$, we fix the element count at $N_x=200$. The error tolerance for all tests is set at $tol=10^{-6}$, which would be well-suited for typical practical applications. Figure \ref{fig:mgscaling} shows the iteration count and computational scaling when either the horizontal order or element count is increased. Both methods shows excellent scalability, with the iteration count showing very little variation, with no tendency to increase at higher DoFs. Moreover, letting $n$ denote the DoF, both methods achieve $O(n)$ computational scaling, independent of whether the increase in DoF is due to element count or polynomial order being increased. 

\section{Conclusion}
The spectral element method has been shown to solve the incompressible Navier-Stokes equations with a free surfaces, achieving spectral convergence of errors at both varying depth and wave steepness. The method provides the opportunity to balance computational efficiency with high numerical accuracy when simulating the temporal evolution of waves over long time periods. The method has also been shown to match experimental results for wave propagation over a bar. Moreover, a geometric $p$-multigrid accelerated iterative solver based on methods such preconditioned defect corrections or a Krylov-based GMRES has been shown to effectively solve the occurring mixed-stage Poisson Pressure problem. The solver is demonstrated to enable fast computations in two-space dimensions and achieves $O(n)$ computational scaling, both when increasing element count and polynomial discretization order. The iterative solver strategy is conceptually straightforward to extend to three space dimensions and is highly suitable for massively parallel implementations, e.g. see related works on FNPF-based free surface models \cite{EngsigKarupEtAl2012,GlimbergEtAl2019, AllanLaskowski}.

In ongoing work, the new free surface INS model is to be extended to wave-structure applications in three space dimensions utilizing a key advantage of a SEM-based solver strategy, namely, the accurate representation of complex geometry via unstructured meshes. 

\bibliographystyle{elsarticle-num-names} 
\bibliography{sn-article}

\appendix
\section{Generation and absorption zones}\label{app:waveZones}
To generate waves through the use generation and absorption zones, a relaxed solution $u^*$ is introduced. The relaxed solution is a combination of two relaxation functions $\Gamma_a(x)$ and $\Gamma_g(x)$ acting on the solution $u(x)$ and an analytical solution $u_e(x)$, 
\begin{align}
    u^*(x) = \Gamma_a(x)u(x)+\Gamma_g(x)u_e(x),
\end{align}
The relaxation functions are piecewise functions, defined as
\begin{subequations}
    \begin{align}
        \Gamma_g (x) & = 
        \begin{cases} 
        f_g(1-y), \quad  & x < 0 \\
        0, \quad &  0 \leq x < 35 \\
        f_g(y), \quad & 35 \leq x
        \end{cases}, \\
        \Gamma_a (x) & = 
        \begin{cases} 
        f_g(y), \quad  & x < 0  \\
        1, \quad &  0 \leq x < 25 \\
        f_a(1-y), \quad &  25 \leq x < 35 \\
        0, \quad &  35 \leq x
        \end{cases},
    \end{align}
\end{subequations}
with relaxation functions defined to guarantee smooth transitions across the relaxation zone interfaces as derived in \cite{ENG06}
\begin{subequations}
    \begin{align}
        f_g(y) &= -2y^3 + 3y^2,  \\
        f_a(y) &= 1 - (1-y)^5. 
    \end{align}
\end{subequations}
Let $x_1,x_2\in \mathbb{R}$ represent the spatial bounds of the domain where $x\in[x_1,x_2]$. Here $y$ is the mapping $y: x \mapsto [0, 1]$.

\end{document}